\documentclass[11pt]{article}
\usepackage[utf8]{inputenc}
\usepackage{amsmath,amssymb,epsfig,bbm}
\usepackage{stmaryrd,mathabx,textcomp}
\usepackage{comment}
\usepackage{color}
\usepackage[T1]{fontenc}

\usepackage{interval}
\usepackage{float}
\usepackage[percent]{overpic}

\usepackage[textsize=small]{todonotes}
\usepackage{enumitem}
\usepackage{varwidth}
\setlist{nolistsep}
\usepackage[colorlinks]{hyperref}

\newtheorem{defi}{Definition}[section]
\newtheorem{prop}[defi]{Proposition}
\newtheorem{theo}[defi]{Theorem}
\newtheorem{theofr}[defi]{Théorème}
\newtheorem{conj}[defi]{Conjecture}
\newtheorem{lemm}[defi]{Lemma}
\newtheorem{lemmfr}[defi]{Lemme}
\newtheorem{coro}[defi]{Corollary}
\newtheorem{rema}[defi]{Remark}
\newtheorem{exem}[defi]{Example}
\newtheorem{exems}[defi]{Examples}

\newcommand{\bdefi}{\begin{defi}}
\newcommand{\edefi}{\end{defi}}
\newcommand{\bprop}{\begin{prop}}
\newcommand{\eprop}{\end{prop}}
\newcommand{\btheo}{\begin{theo}}
\newcommand{\etheo}{\end{theo}}
\newcommand{\btheofr}{\begin{theofr}}
\newcommand{\etheofr}{\end{theofr}}
\newcommand{\blemm}{\begin{lemm}}
\newcommand{\elemm}{\end{lemm}}
\newcommand{\blemmfr}{\begin{lemmfr}}
\newcommand{\elemmfr}{\end{lemmfr}}
\newcommand{\brema}{\begin{rema}}
\newcommand{\erema}{\end{rema}}
\newcommand{\bexer}{\begin{exem}}
\newcommand{\eexer}{\end{exem}}
\newcommand{\bexems}{\begin{exems}}
\newcommand{\eexems}{\end{exems}}
\newcommand{\bconj}{\begin{conj}}
\newcommand{\econj}{\end{conj}}
\newcommand{\bcoro}{\begin{coro}}
\newcommand{\ecoro}{\end{coro}}
\newcommand{\dem}{\noindent{\bf Proof. }}

\usepackage{mathrsfs}
\renewcommand\mathcal{\mathscr}

\newcommand{\A}{{\cal A}}

\newcommand{\F}{{\cal F}}

\newcommand{\N}{{\cal N}}
\newcommand{\OOO}{{\cal O}}

\newcommand{\R}{{\cal R}}
\newcommand{\maths}[1]{{\mathbb #1}}  

\newcommand{\CC}{\maths{C}}

\newcommand{\NN}{\maths{N}}

\newcommand{\QQ}{\maths{Q}}
\newcommand{\RR}{\maths{R}}

\newcommand{\ZZ}{\maths{Z}}

\newcommand{\ra}{\rightarrow}

\newcommand{\ov}[1]{{\overline #1}} 
\newcommand{\wt}[1]{{\widetilde{#1}}}
\newcommand{\wh}[1]{{\widehat{#1}}}

\newcommand{\ga}{\gamma}

\newcommand{\cqfd}{\hfill$\Box$}

\newcommand{\bigO}{\operatorname{O}}

\newcommand{\card}{{\operatorname{Card}}}

\newcommand{\covol}{\operatorname{covol}}

\newcommand{\diam}{{\operatorname{diam}}}

\newcommand{\Leb}{\operatorname{Leb}}

\newcommand{\sg}{\operatorname{sg}}

\newcommand{\ssm}{\!\smallsetminus\!}

\newcommand{\tr}{\operatorname{\tt tr}}
\newcommand{\n}{\operatorname{\tt n}}

\newcounter{fig}

\title{On the multiplicative pair correlations\\ of sums of two squares}
\author{Jouni Parkkonen \and Fr\'ed\'eric Paulin}
\date{\today}

\begin{document}
\bibliographystyle{../alphanum}
\maketitle
\begin{abstract} 
We study the pair correlations of the logarithms of the integral
values of quadratic norm forms at various scalings, proving the
existence of pair correlation measures. We describe a surprising set
of asymptotic behaviours when the scaling increases, passing from a
punctual measure to a Poissonian behaviour through an exotic behaviour
at the transition phase.
\end{abstract}

\section{Introduction}
\label{sec:intro}

In this paper, we study the asymptotic distribution of the
(generalized) pair correlations of the logarithms of the values on
integer points of integral quadratic norm forms, including the special
case of sums of two squares.  See for instance \cite{Hooley94,
  FreKurRos17,BerMoo25} and their references for various
distributional aspects of sums of two squares.

A framework for this study, that includes the one in \cite{NaiPol07,
  HauZaf24}, is the following one, see also \cite{ParPau22MJCNT,
  ParPau24RNT,Sayous25, Sayous26}. Let $\F=(F_N,\omega_N)_{N\in\NN}$
be a nondecreasing sequence of finite subsets $F_N$ of $\RR$, endowed
with a {\it weight function} $\omega_N: F_N \ra\; ]\,0,+\infty\,[\,$.
Let $\phi: \NN \ra\; ]\,0,+\infty\,[$ be a {\it scaling function}
converging to $+\infty$.  Let $\psi: \NN \ra\; ]\,0,+\infty\,[$ be a
{\it renormalizing function}, that will be naturally chosen depending
on $\phi$. We denote by $\Delta_z$ the unit Dirac measure at any point
$z$ of any measurable space. We define the {\it empirical pair
  correlation measure of $\F$ at time $N$ with scaling $\phi(N)$} as
the measure on $\RR$ with finite support
\[
\R^{\F,\phi}_N=\frac{1}{\psi(N)}\sum_{x,y\in F_N\;:\;x\neq y}\;
\omega_N(x)\,\omega_N(y)\,\Delta_{\phi(N)(y-x)}\,.
\]
When the sequence of measures $(\R^{\F,\phi}_N)_{N\in\NN}$ weak-star
converges to a measure $m_{\F,\phi}$ on $\RR$, we call $m_{\F,\phi}$
the {\it asymptotic pair correlation measure of $\F$ for the scaling}
$\phi$. When $m_{\F,\phi}=\rho_{\F,\phi}\,\Leb_\RR$ is absolutely
continuous with respect to the Lebesgue measure $\Leb_\RR$ of $\RR$,
the Radon-Nikodym derivative $\rho_{\F,\phi}$ is called the {\it
  asymptotic pair correlation density} of $\F$ {\em for the scaling}
$\phi$. When $\rho_{\F,\phi}$ is constant, ones says that the pair
correlations exhibit a {\it Poissonian asymptotic behaviour}.

Let $K$ be a quadratic imaginary number field, with discriminant
$D_K$, ring of integers $\OOO_K$ and (relative) norm $\n:z\mapsto
z\,\overline{z}$. For every $a\in\NN$, let $r_K(a)=\card\{z\in\OOO_K:
\n(z)=a\}$ be the number of representations of $a$ by the norms of
elements of $\OOO_K$. We fix $\alpha\in\;]0, \frac{1}{2}[$ throughout
this paper. For all nonzero $a,b,N\in\NN$, let
\begin{equation}\label{eq:defirKalpha}
r_{K,\alpha}(a,b,N)=\card\{(w,z)\in\OOO_K:
\n(w)=a, \;\n(z)=b,\; \n(z-w)\leq N^{2\alpha}\}\,,
\end{equation}
which if $N$ is large enough is equal to the product $r_K(a)\,r_K(b)$
of the numbers of representations by the norm form $\n$ of $a$ and of
$b$. In this paper, we study the asymptotic behaviour of the following
empirical distribution of pair correlations
\begin{equation}\label{eq:defiRsubN}
\R_{N}=\frac{1}{\psi(N)}
\sum_{a,b\in \NN\;:\;  a\neq b,\;0<a,b\le N^2}\;
r_{K,\alpha}(a,b,N)\,\Delta_{\phi(N)(\ln a-\ln b)}\,.
\end{equation}
In particular, when $K=\QQ(i)$, $D_K=-4$, the behaviour of the
measures $\R_{N}$ as $N\ra+\infty$ gives the {\it multiplicative pair
  correlation} asymptotics of the sums of two squares, more precisely,
the asymptotic of the values $(\frac{a}{b})^{\phi(N)}$ of the ratios
$\frac{a}{b}$ of the sums of two squares $a$ and $b$ raised to the
power $\phi(N)$, weighted by the number of their representations in
sectors.

In addition to the purely arithmetic interest, a geometric motivation
in order to study these representations is that the logarithms of the
norms of $\OOO_K$, when $K$ has class number one, form the ortholength
spectrum of geodesic segments from a neighbourhood of the (unique)
cusp to itself in the $3$-dimensional real hyperbolic Bianchi orbifold
of $K$, see \cite[Sect.~7]{ParPau24RNT} for details.

In order to simplify the statements in this introduction, we only
consider the power scalings $\phi:N\mapsto N^\beta$ for $\beta\in\;
]0,1+\frac{\alpha}{2}[\,$, where this upper bound on $\beta$
corresponds to the restriction in Equation \eqref{eq:defirKalpha} on
the pairs of representations. Such a power scaling is a usual choice,
see for instance \cite{NaiPol07,HauZaf24}.  We denote by
$\R_N^{\alpha,\beta}$ the empirical distribution given by Equation
\eqref{eq:defiRsubN} for this power scaling. We define a continuous,
positive, piecewise real analytic, even function $\rho_{1-\alpha}:
\RR\ra[0,+\infty[$ by
{\renewcommand{\arraystretch}{1.4}
%%%
\begin{equation}\label{eq:defpaircorrfunct}
t\mapsto\left\{\begin{array}{ll}
\frac{8\,\pi}{3\,|D_K|} &\text{if } t=0
\\
\frac{8\,\pi}{|D_K|\,t^3}\Big(\arcsin(\frac{t}{2})-
\frac{t}{2}(1-\frac{t^2}{2})\sqrt{1-\frac{t^2}{4}}\;\Big)
&\text{if } 0<|t|\leq 2 \\
\frac{4\,\pi^2}{|D_K|\,|t|^3}
&\text{if } |t|> 2\,.
\end{array}\right.
\end{equation}
}

\btheo\label{theo:mainintro} Assume that we have $D_K\equiv 0\bmod
4$. Let $\alpha\in\;]0, \frac{1}{2}[$ and $\beta\in\;]0,1+
\frac{\alpha}{2}[\,$.  As $N\ra+\infty$, the empirical pair
correlation measures $\R_N^{\alpha,\beta}$ converge, for the weak-star
convergence of measures on $\RR$, to the asymptotic pair correlation
measure $m_{\alpha,\beta}$ given by
\[
m_{\alpha,\beta}=\begin{cases}
\frac{4\,\pi^2}{|D_K|}\;\Delta_0
&\textrm{if } \beta=\;]0,1-\alpha[ \text{ and } \psi(N)=N^{2+2\alpha},\\
\rho_{1-\alpha}\;\Leb_{\RR}
&\textrm{if } \beta=1-\alpha \text{ and }
\psi(N)=N^{2+2\alpha}=N^{3+\alpha-\beta},\\
\frac{8\pi}{3\,|D_K|}\;\Leb_{\RR}
&\textrm{if } \beta\in\;]1-\alpha,1+\frac{\alpha}{2}[\;,\;
  \alpha\leq \frac{1}{6} \text{ and } \psi(N)=N^{3+\alpha-\beta}.
\end{cases}
\]
\etheo

In particular, the asymptotic behaviour of the empirical pair
correlation measures $\R_N^{\alpha,\beta}$ is Poissonian when
$\beta\in\; ]1-\alpha,1+ \frac{\alpha}{2}[$. It has a phase transition
with an exotic asymptotic pair correlation density $\rho_{1-\alpha}$
with respect to the Lebesgue measure when $\beta=1-\alpha$, see
Figure \hyperlink{figintro}{\addtocounter{fig}{1}\arabic{fig}%
\addtocounter{fig}{-1}} below for an example. Below this
threshold, the empirical pair correlation measures concentrate on a
punctual measure.  Such phase transition phenomenona as the scaling
increases appear frequently, see for instance \cite{ParPau22MJCNT,
ParPau24RNT, Sayous25,Sayous26}.

We refer to Theorem \ref{theo:sumsquaremain} for a more complete
version of Theorem \ref{theo:mainintro}, without the restriction on
the discriminant $D_K$, with more general scaling functions, as well
as for error terms. These error terms constitute the main technical
parts of this paper. It is an interesting feature that even in the
Poissonian asymptotic behaviour case (when $\beta\in\;]1-\alpha,
1+\frac{\alpha}{2}[\,$), the validity of the error terms depend on
whether $\beta\in\;]1-\alpha, 1-\frac{\alpha}{2}]$, $\beta\in
]1-\frac{\alpha}{2}, 1[$ or $\beta\in[1,1+\frac{\alpha}{2}[\,$, while
the constant value $\frac{8\pi}{3\,|D_K|}$ of the asymptotic pair
correlation density does not change.

\begin{center}
\hypertarget{figintro}{\includegraphics[width=14cm]{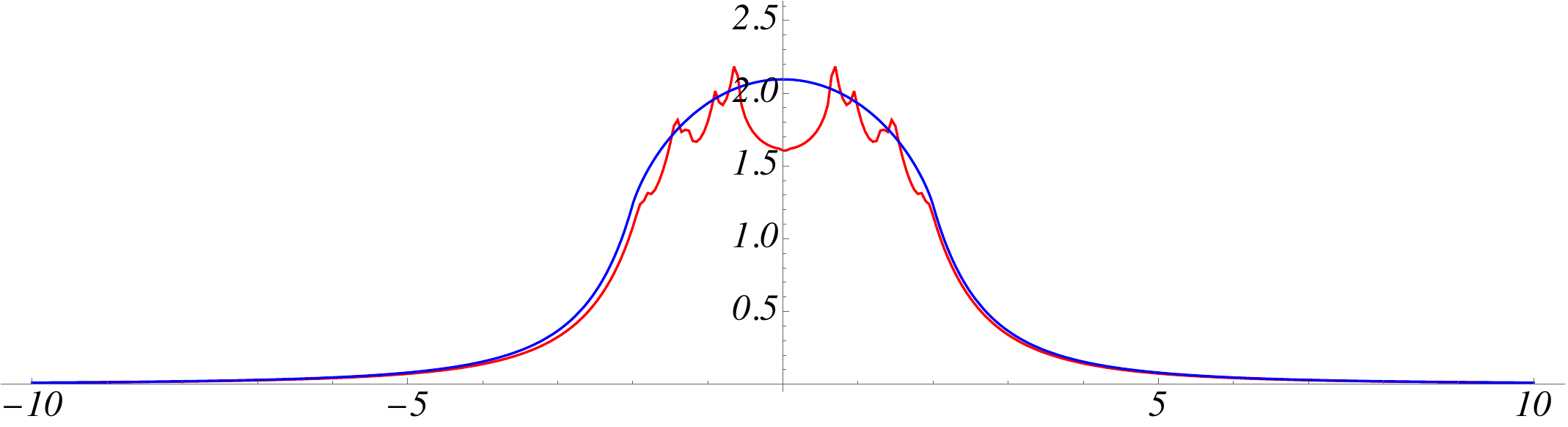}
\includegraphics[width=14cm]{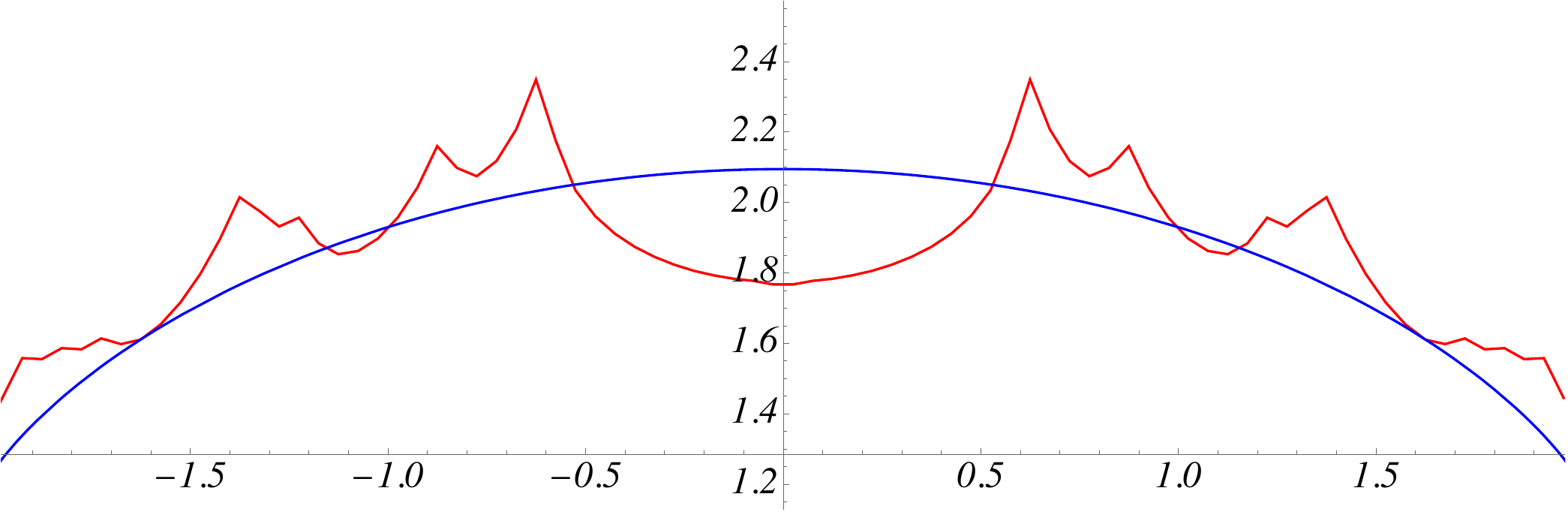}}
\end{center}\vspace*{-0.5cm}
\begin{center}
\addtocounter{fig}{1} {\small Figure \arabic{fig}~: For $\alpha=0.15$,
  $\beta=0.85$, graph of $\rho_{1-\alpha}$ (in blue) with the density
  of $\R_{N}^{\alpha,\beta}$ (in red) on $[-10,10]$ for $N=2000$,
  and on the shorter interval $[-2,2]$ for $N=3000$.}
\end{center}

\medskip
The study of pair correlations in a noncompact setting has a rich
history, including the seminal paper \cite{Montgomery73} on the zeros
of the Riemann zeta function.  The lengths of the closed geodesics in
negative curvature have a Poissonian pair correlation asymptotics or
their empirical distribution converges to an exponential probability
measure, depending on the scaling factor, see \cite{PolSha06,
  ParPau23BSMF}.  For real numbers $\alpha',\beta'$ satisfying some
Diophantine condition, the image of $\ZZ^2$ by the quadratic form
$(x,y)\mapsto (x-\alpha')^2+(y-\beta')^2$ also exhibits a Poissonian
pair correlation asymptotic behaviour, by \cite{Marklof03}.  Similar
problems often arise in quantum chaos, including energy level spacings
or clusterings, and in statistical physics, including molecular
repulsion or interstitial distribution, and in various
number-theoretical contexts. See for instance \cite{Berry88, RudSar98,
  Vanderkam99, Zaharescu03, ElkMcM04, BocZah05, MarStr13, ElBMarVin15,
  LarSto20, HofKal21, Weiss23, LutSouTec25}.

In Section \ref{rem:casBetageq2}, we discuss the results of
experiments that indicate a change from the Poissonian asymptotic
behaviour of the pair correlations when $\beta$ is beyond the range of
Theorem \ref{theo:mainintro}, with strong level repulsion phenomena,
and finally a total loss of mass for $\beta>2$.  It would also be
interesting to study the weighted family
\[
\F_K=\big(F_N=\{\ln n\;:\; 0<n \leq N, \;r_{K}(n)\neq 0\},
w_N=r_{K}\circ \exp\big)_{N\in\NN}
\]
with weights given by $r_{K}$ (removing the zero weights). For
instance, when $D_K=-4$, the function $r_K$ is the number of
representations of integers as sums of two squares. In
\cite[Coro.~2.5]{ParPau24RNT}, we proved that the asymptotic pair
correlation density of $\F_K$ with constant scaling function $\phi=1$,
i.~e.~the distribution function of $\F_K$, is $t\mapsto \frac 12
e^{-|t|}$.  See also \cite{ParPau23BSMF} for a general result when
$\phi$ is constant.

\medskip
\noindent{\small {\it Acknowledgements: } This research was supported
  by the French-Finnish CNRS IEA PaCap.}

\section{Parametrizing the pair correlation data}
\label{sec:unfold}

Let us fix $N\in\NN\ssm\{0\}$. In this section, we give a
parametrisation of the set of representations $q\in\OOO_K$ by the norm
form $\n$ of the positive integers $a$ at most $N$, under the
additional constraints for them to be part of the pair correlation
data appearing in the empirical distribution $\R_N$ given by Equation
\eqref{eq:defiRsubN}.

We consider the (relative) {\it trace} and {\it norm} maps from $\CC$
to $\RR$ defined by $\tr:z\mapsto z+\ov{z}$ and $\n:z\mapsto
z\,\ov{z}$ respectively. Let $(1,\omega_K)$ be the usual $\ZZ$-basis
of $\OOO_K$, with
\[
\omega_K = \left\{\begin{array}{ll} \!\!\!\frac{i\sqrt{|D_K|}}{2}&
\!\text{if~} D_K\equiv 0\!\!\mod 4\\ \!\!\!\frac{1+i\sqrt{|D_K|}}{2}&
\!\text{otherwise.}\end{array}\right.
\]
We have 
{\renewcommand{\arraystretch}{1.3}
\[
\tr(\omega_K)=
\left\{\begin{array}{ll} \!\!0& \!\text{if~} D_K\equiv 0\!\!\mod
4\!\\ \!\!1& \!\text{otherwise,}\end{array}\right.
\quad {\rm and }\quad \n(\omega_K)=
\left\{\begin{array}{ll}\!\!\frac{|D_K|}{4}& \!\text{if~} D_K\equiv
0\!\!\mod 4\!\\ \!\!\frac{1+|D_K|}{4}& \!\text{otherwise.}\end{array}
\right.
\]
}%
The norm of any nonzero element of $\OOO_K$ is a positive integer. In
particular, the area $\covol_{\OOO_K}$ of the fundamental
parallelogram $[0,1]+ \omega_K[0,1]$ of the $\ZZ$-lattice $\OOO_K$ in
$\CC$, its diameter $\diam_{\OOO_K}$ and the shortest length
$\operatorname{Sys}_{\OOO_K}$ of a nonzero element of $\OOO_K$ satisfy
\begin{equation}\label{eq:covoldiamsys}
  \covol_{\OOO_K}=\frac{\sqrt{|D_K|}}{2},\qquad
  1\leq \diam_{\OOO_K}=\bigO(\sqrt{|D_K|})\quad\text{and}\quad
  \operatorname{Sys}_{\OOO_K}=1\,.
\end{equation} 

In this section, we fix $p\in \OOO_K\ssm\{0\}$ and we write
\[
p=x_p+\omega_K\,y_p
\]
with $x_p,y_p\in\ZZ$. We define
\[
x'_p= 2\,x_p+y_p\,\tr(\omega_K) \qquad\text{and}\qquad
y'_p= x_p\,\tr(\omega_K)+2\,\n(\omega_K)\,y_p\,,
\]
which are easily seen to be elements of $\ZZ$ that are not
simultaneously zero. We also denote by $c'_p=(x'_p,y'_p)\in\NN\ssm
\{0\}$ the (positive) greatest common divisor of $x'_p$ and $y'_p$,
and we define
\[
v_p=\frac{1}{c'_p}(y'_p-\omega_K\,x'_p)\in\OOO_K\ssm\{0\}\,.
\] 
Note that when $D_K\equiv 0\!\!\mod 4$, then $\tr(\omega_K)=0$ and
$\overline{\omega_K}=-\omega_K$, hence $x'_p= 2\,x_p$, $y'_p=
2\,\n(\omega_K)\,y_p$ and
\begin{equation}\label{eq:valsppvpdiscnul}
  c'_p\,v_p=2\,\n(\omega_K)\,y_p-2\,x_p\,\omega_K=
  -2\,\omega_K(x_p+\omega_K\,y_p)=
  -i\,\sqrt{|D_K|}\;p\,.
\end{equation}
Note that there exists a constant $c_K>0$ depending only on $K$ such
that
\begin{equation}\label{eq:controlcppvp}
  1\leq \max\{c'_p,|v_p|\}\leq c_K\,|p|\qquad\text{and}\qquad
  \frac{1}{c_K}\,|p|\leq c'_p\,|v_p|\leq c_K\,|p|\,.
\end{equation}

We denote by $H^+_{p}$ the open halfplane containing $0$ whose
boundary is the mediatrix $-\frac{p}{2}+i\,p\,\RR$ of the segment
$[0,-p]$, by $B(-p,N)$ the closed ball of center $-p$ and radius $N$,
and we define
\[
J_{p,N}= \big(\OOO_K\cap H^+_{p}\cap B(-p,N)\big)\ssm\{0\}\,.
\]
We also define
\[
\Lambda_p=\{\tr(\ov{p}\,q)\;:\;q\in \OOO_K\}\subset\ZZ\,,
\]
which is a $\ZZ$-lattice of $\RR$ by the linearity (over $\RR$) of the
trace. For every $t\in\RR$, let
\[
L_{p,t}= \{z\in\CC:\tr(\ov{p}\,z)=t\,c'_p\}\,,
\]
which is an affine (real) line of $\CC$, again by the linearity of the
trace. The set of blue dots below represents $J_{p,N}$ with $\OOO_K=
\ZZ[i]$, $p=2+3\,i$, $c'_p=2$, $v_p= 3-2i$, and the red lines are
$L_{p,k}$ with $k=0$ and $k=19$.

\hypertarget{figure2}{

\addtocounter{fig}{1}
\begin{center}
\begin{picture}(0,0)%
\includegraphics{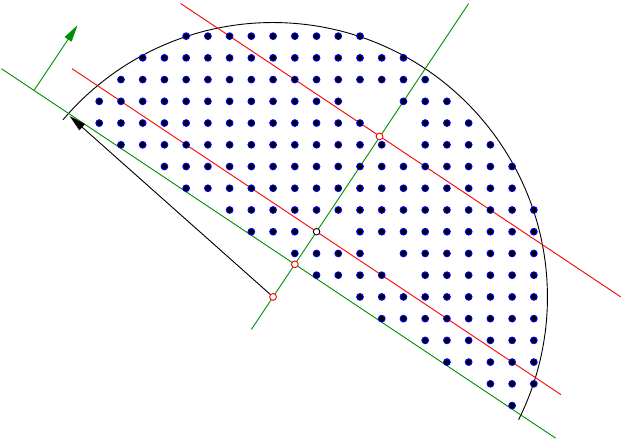}%
\end{picture}%
\setlength{\unitlength}{3812sp}%
\begin{picture}(5154,3624)(-641,-2773)
\put(4141,-1321){\makebox(0,0)[lb]{\smash{\fontsize{11}{13.2}\usefont{T1}{ptm}{m}{n}{\color[rgb]{0,0,0}$L_{p,k}$}%
}}}
\put(1600,-1714){\makebox(0,0)[lb]{\smash{\fontsize{11}{13.2}\usefont{T1}{ptm}{m}{n}{\color[rgb]{0,0,0}$-p$}%
}}}
\put(2353,-583){\makebox(0,0)[lb]{\smash{\fontsize{11}{13.2}\usefont{T1}{ptm}{m}{n}{\color[rgb]{0,0,0}$p$}%
}}}
\put(4051,-2401){\makebox(0,0)[lb]{\smash{\fontsize{11}{13.2}\usefont{T1}{ptm}{m}{n}{\color[rgb]{0,0,0}$L_{p,0}$}%
}}}
\put(-494,479){\makebox(0,0)[lb]{\smash{\fontsize{11}{13.2}\usefont{T1}{ptm}{m}{n}{\color[rgb]{0,0,0}$H^+_p$}%
}}}
\put(  1,-511){\makebox(0,0)[lb]{\smash{\fontsize{11}{13.2}\usefont{T1}{ptm}{m}{n}{\color[rgb]{0,0,0}$N$}%
}}}
\put(1486,-1411){\makebox(0,0)[lb]{\smash{\fontsize{11}{13.2}\usefont{T1}{ptm}{m}{n}{\color[rgb]{0,0,0}$-\frac{p}{2}$}%
}}}
\put(3196,569){\makebox(0,0)[lb]{\smash{\fontsize{11}{13.2}\usefont{T1}{ptm}{m}{n}{\color[rgb]{0,0,0}$\RR\,p$}%
}}}
\put(2251,-72){\makebox(0,0)[lb]{\smash{\fontsize{11}{13.2}\usefont{T1}{ptm}{m}{n}{\color[rgb]{0,0,0}$w_{p,k}$}%
}}}
\put(2551,-259){\makebox(0,0)[lb]{\smash{\fontsize{11}{13.2}\usefont{T1}{ptm}{m}{n}{\color[rgb]{0,0,0}$z_{p,k}$}%
}}}
\put(2509,-1360){\makebox(0,0)[lb]{\smash{\fontsize{11}{13.2}\usefont{T1}{ptm}{m}{n}{\color[rgb]{0,0,0}$v_p$}%
}}}
\put(2038,-1075){\makebox(0,0)[lb]{\smash{\fontsize{11}{13.2}\usefont{T1}{ptm}{m}{n}{\color[rgb]{0,0,0}$0$}%
}}}
\end{picture}%

\end{center}
\vspace*{-0.8cm}
\begin{center}
{\small Figure \arabic{fig}~: A description of the finite set $J_{p,N}$.}
\end{center}
}

Recall that we have 
\begin{equation}\label{eq:tracescalarproduct}
  (z=x+i\,y,z'=x'+i\,y')\mapsto \tr(\overline{z}\,z')=
  2(x\,x'+y\,y')=2\,\langle z,z'\rangle\,,
\end{equation}
where $\langle\;,\; \rangle$ is the usual scalar product of the
Euclidean real plane $\CC$. The following result summarizes the
geometric properties of the above objects.

\blemm\label{lem:calculavecp} Let $N\in\NN \ssm\{0\}$ and $p\in\OOO_K
\ssm \{0\}$.
\begin{enumerate}
\item\label{item1:calculavecp} The nonzero vectors $v_p$ and $p$ are
  perpendicular in $\CC$~: we have $\tr(\ov{p}\,v_p)=0$ and $v_p\in
  L_{p,0}$.
\item\label{item2:calculavecp} We have
  $H^+_{p}=\{z\in\CC\;:\;\n(z)<\n(p+z)\}$, hence \[J_{p,N}=
  \{q\in\OOO_K\;:\; 0< \n(q)<\n(p+q)\leq N^2\}\,.\]
  Furthermore, $J_{p,N}$ is empty if $|p|\geq 2N$.
\item\label{item3:calculavecp}
  We have $\Lambda_p=c'_p\,\ZZ$.
\item\label{item4:calculavecp} For every $t\in\RR$, the affine (real)
  line $L_{p,t}$ is perpendicular to the (real) line $\RR\, p$ (hence
  is parallel to the boundary $\partial H^+_{p}$ of $H^+_{p}$) and
  meets $\RR\, p$ exactly at the point $z_{p,t}=\frac{t\,c'_p}{2\n(p)}
  \,p$ (see Figure \hyperlink{figure2}{\arabic{fig}}).
  \item\label{item5:calculavecp} For every $t\in\RR$, the intersection
    $L_{p,t}\cap\OOO_K$ is nonempty if and only if $t\in\ZZ$. For
    every $k\in\ZZ$, if $w_{p,k}$ is one of the at most two points of
    $L_{p,k}\cap \OOO_K$ the closest to $z_{p,k}$ (see Figure
    \hyperlink{figure2}{\arabic{fig}}), we have $L_{p,k}\cap \OOO_K=
    w_{p,k}+\ZZ\,v_p$. Furthermore, there exists $t_{p,k}\in
    [-\frac12,\frac12]$ such that $w_{p,k} -z_{p,k}=t_{p,k}\,v_p$.
\item\label{item6:calculavecp} Let $\kappa_{p} = \big\lfloor
  -\frac{\n(p)}{c'_p}\rfloor+1$. For every $k\in\ZZ$, the affine line
  $L_{p,k}$ is contained in the open halfspace $H^+_p$ if and only if
  we have $k\geq \kappa_{p}$.
\item\label{item7:calculavecp} The set $J_{p,N}$ is the set of
  elements $\frac{k\,c'_p}{2\,\n(p)}\,p+(t_{p,k}+\ell)v_p$ such that
  $k,\ell\in\ZZ$, $k\geq \kappa_{p}$ and $\big(
  \frac{k\,c'_p}{2\,\n(p)}+1\big)^2|p|^2+(t_{p,k}+\ell)^2|v_p|^2 \leq
  N^2$.
\end{enumerate}
\elemm

\dem For every $z=x+\omega_Ky\in \CC$ with $x,y\in\RR$, by the
definition of $x'_p$ and $y'_p$, we have
\[
\tr(\ov{p}\,z)= 2(x_p\,x+\n(\omega_K)\,y_p\, y)+
(x_p\,y+y_p\,x)\,\tr\omega_K =x'_p\,x+y'_p\,y\,.
\]
This proves Assertion \eqref{item1:calculavecp} by the definition of
$v_p$ and by Equation \eqref{eq:tracescalarproduct}.  This also proves
Assertion \eqref{item3:calculavecp} since then $\Lambda_p=x'_p\,\ZZ+
y'_p\,\ZZ =c'_p\,\ZZ$ by the definition of $c'_p$.

The set $\{z\in\CC:\langle p,z+\frac{p}{2}\rangle>0\}$ is the
halfplane in $\CC$ containing $0$ with boundary the affine line
$-\frac{p}{2}+ p^\perp$, hence is equal to $H^+_{p}$. The first claim
of Assertion \eqref{item2:calculavecp} is then immediate using
Equation \eqref{eq:tracescalarproduct} since for every $z\in\CC$, we
have $\n(z)<\n(p+z)$ if and only if $\tr(\ov{p}\, (z+\frac{p}{2}))
>0$. The second claim follows since $B(-p,N)= \{z\in\CC:|z+p|\leq N\}$
and the norm $\n$ is the square of the absolute value.  The
intersection $H^+_{p} \cap B(-p,N)$ is empty if $|p|\geq 2N$, as the
inequalities $\n(z)<\n(p+z)\leq N^2$ imply that $|z|<N$ and that
$|z|\geq |p|-N$ by the inverse triangle inequality. The last claim of
Assertion \eqref{item2:calculavecp} follows.

For all $\lambda,t\in\RR$, we have $\lambda\, p\in L_{p,t}$ if and
only if $\tr(\ov{p}\,\lambda\, p)=t\,c'_p$, that is, if and only if
$\lambda=\frac{t\,c'_p}{2\n(p)}$. Assertion \eqref{item4:calculavecp}
then follows, since $L_{p,0}=\{z\in\CC:\tr(\ov{p}\,z\}=0\}$ is the
vector line orthogonal to $\RR\,p$ and $L_{p,t}=z_{p,t}+L_{p,0}$.

Let us prove Assertion \eqref{item5:calculavecp}. By the definitions
of $\Lambda_p$ and $L_{p,t}$, the intersection $L_{p,t}\cap\OOO_K$ is
nonempty if and only if $t\,c'_p\in \Lambda_p$, hence the first claim
follows from Assertion \eqref{item3:calculavecp}. Let $k\in\ZZ$.
Since the point $w_{p,k}$ belongs by construction to the affine line
$L_{p,k}$, which is perpendicular to $\RR\, p$ by Assertion
\eqref{item4:calculavecp}, and since the vector $v_p$ is nonzero and
perpendicular to $p$ by Assertion \eqref{item1:calculavecp}, we have
$L_{p,k}=w_{p,k} +\RR\,v_p$.  Since $w_{p,k}\in \OOO_K$, we have
$L_{p,k}\cap\OOO_K= \{w_{p,k}+s\,v_p : s\,v_p\in \OOO_K\}$. The second
claim of Assertion \eqref{item5:calculavecp} then follows from the
fact that the vector $v_p\in\OOO_K\ssm \{0\}$ is primitive (it has
relatively prime integral coordinates in the $\ZZ$-basis $(1,
\omega_K)$ of $\OOO_K$). This proves that two consecutive points of
$\OOO_K$ on the affine line $L_{p,k}$ are at distance exactly $|v_p|$.
Hence we have $|w_{p,k} -z_{p,k}|\leq \frac{1}{2}\,|v_p|$ for every
$k\in\ZZ$, and the last claim of Assertion \eqref{item5:calculavecp}
follows.

Let us now prove Assertion \eqref{item6:calculavecp}. The affine line
$\partial H^+_p$, which contains $-\frac{p}{2}$ and is perpendicular
to $\RR\,p$, is equal to 
\[
\Big\{z\in \CC:\tr(\ov{p}\,z)= \tr\Big(\ov{p}\,
\Big(\!-\frac{p}{2}\,\Big)\Big) =-\n(p)\Big\}=L_{p,-\frac{\n(p)}{c'_p}}\,.
\]
Therefore for every $k\in\ZZ$, the affine line $L_{p,k}$ is contained
in $H^+_p$ if and only if $k>-\frac{\n(p)}{c'_p}$, that is, if and
only if $k\geq \kappa_{p}= \big\lfloor -\frac{\n(p)}{c'_p}\,
\big\rfloor+1$. Note that $\kappa_{p}$ is nonpositive.

Finally, let us prove Assertion \eqref{item7:calculavecp}.  Note that
$\CC$ is foliated by the affine lines $L_{p,t}$ for $t\in\RR$.  For
every $t\in\RR$, since $J_{p,N}$ is contained in $\OOO_K$, the
intersection $L_{p,t} \cap J_{p,N}$ is empty if $t\notin\ZZ$ by
Assertion \eqref{item5:calculavecp}. Hence $J_{p,N}=\bigcup_{k\in\ZZ}
L_{p,k} \cap J_{p,N}$. By Assertion \eqref{item6:calculavecp} and the
definition of $J_{p,N}$, we hence have $J_{p,N}=\bigcup_{k\geq
  \kappa_p} L_{p,k} \cap \OOO_K \cap B(-p,N)$. By Assertion
\eqref{item5:calculavecp}, any element $z\in L_{p,k} \cap \OOO_K$ can
be uniquely written as $z=w_{p,k}+\ell\,v_p$ for some $\ell\in\ZZ$.
Hence by Assertions \eqref{item4:calculavecp} and
\eqref{item5:calculavecp}, we have
\[
z=z_{p,k} +w_{p,k}-z_{p,k}+\ell\,v_p=
\frac{k\,c'_p}{2\n(p)}\, p+t_{p,k}\,v_p+\ell\,v_p\,.
\]
This proves the result, since $p$ and $v_p$ are orthogonal, hence the
inequality $|z+p|\leq N$ is equivalent to the last inequality of
Assertion \eqref{item7:calculavecp}. \cqfd

\bigskip
Let us consider the map $j_{p,N}:\RR^2\ra\CC$ defined by
\begin{equation}\label{eq:defijpN}
  j_{p,N}:(s,t)\mapsto N\Big(s\,\frac{p}{|p|}+t\,\frac{v_p}{|v_p|}\Big)\,,
\end{equation}
which is a homothety of Euclidean vector spaces (and in particular a
homeomorphism). By Lemma \ref{lem:calculavecp}
\eqref{item7:calculavecp}, we have
\begin{align}
j_{p,N}^{-1}(J_{p,N})=\Big\{&\Big(\,s=\frac{c'_p}{2\,|p|\,N}k,\;\;
t=\frac{t_{p,k}\,|v_p|}{N}+\frac{|v_p|}{N}\ell\,\Big):
\nonumber\\ & k,\ell\in\ZZ,\;\; s> -\frac{|p|}{2\,N},\;\;
\Big(s+\frac{|p|}{N}\Big)^2+t^2\leq 1\Big\}\,.\label{eq:valjpNJpN}
\end{align}
Recall that $\alpha\in\;]0,\frac12[$ has been fixed in the
introduction. The finite subset $j_{p,N}^{-1}(J_{p,N})$ of $\RR^2$ is
contained in $B(-\frac{|p|}{N},1)$ and converges as $N$ tends to
$+\infty$  for the Hausdorff distance on the set of closed subsets of
the metric space $\RR^2$ to the closed halfdisc
\begin{equation}\label{eq:Bplus01}
B^+(0,1)=\{(s,t)\in\RR^2:s\geq 0, \;\;s^2+t^2\leq 1\}\,,
\end{equation}
uniformly in $p\in\OOO_K$ with $0<|p|\leq N^\alpha$ since $\alpha<1$.
Furthermore, since the horizontal coordinate $s$ in Equation
\eqref{eq:valjpNJpN} varies by constant steps $\frac{c'_p}{2\,|p|\,N}$
as $k$ varies in $\ZZ$ and, when $s$ is fixed, the vertical coordinate
$t$ varies by constant steps $\frac{|v_p|}{N}$ as $\ell$ varies in
$\ZZ$, a two-dimensional Riemann sum argument proves that the measure
\[
\frac{c'_p}{2\,|p|\,N}\;\frac{|v_p|}{N}
\sum_{(s,t)\in j_{p,N}^{-1}(J_{p,N})} \Delta_{(s,t)}
\]
weak-star converges as $N\ra+\infty$ (uniformly in $p\in\OOO_K$ with
$0<|p|\leq N^\alpha$) to the restriction ${\rm Leb} _{B^+(0,1)}$ to
$B^+(0,1)$ of the Lebesgue measure of $\RR^2$. In particular, since
the area of $B^+(0,1)$ is $\frac{\pi}{2}$ and by the right-hand part
of Equation \eqref{eq:controlcppvp} for the last equality, as
$N\ra+\infty$ and uniformly in $p\in\OOO_K$ with $0<|p|\leq N^\alpha$,
we have
\begin{equation}\label{eq:cardJpN}
  \card \;J_{p,N} = \card\; j_{p,N}^{-1}(J_{p,N})\sim
  \frac{\pi\,|p|\,N^2}{c'_p\,|v_p|}=\bigO(N^2)\,.
\end{equation}

\section{Uniformisation of the empirical distribution}
\label{sec:linear}

Let us fix again $N\in\NN\ssm\{0\}$, that we will assume to tend to
$+\infty$ after Equation \eqref{eq:netoynupN}. In this technical
section, we represent the empirical pair correlation distribution
$\R_N$ defined in Equation \eqref{eq:defiRsubN} on $[0,+\infty[$ as a
distribution with an explicit density with respect to the Lebesgue
measure, up to a controlled error term, see Equation
\eqref{eq:caspasplusinfty}. We fix in this section
\begin{equation}\label{eq:conditiongaeps}
  \alpha\in\Big]0,\frac{1}{2}\Big[\;,\qquad
      \ga\in\;\Big]0,\frac{1-2\alpha}{2}\Big[\qquad\text{and}\qquad
\epsilon=\epsilon_N=\frac{1}{N^\ga} \in\;]0,1[\,,
\end{equation}
so that $\ga$ exists and $\epsilon\ra0$ as $N\ra+\infty$. We assume in
this section that
\begin{equation}\label{eq:scalingrange}
\lim_{N\ra+\infty}\;\frac{\phi(N)}{N^{2-2\alpha-2\ga}}=0\,.
\end{equation}
When $\phi:N\mapsto N^\beta$ is a power scaling as in the introduction
with $\beta\in\;]0,2-2\alpha[\,$, the above assumption
\eqref{eq:scalingrange} is equivalent to the fact that $\ga\in\;]0,
\frac{2-2\alpha-\beta}{2}[\,$.

\medskip
Let $\log : \CC^\times \ra\CC/(2\pi i \ZZ)$ be the biholomorphic
(complex) logarithm map. Note that the trace map $\tr$, being constant
modulo $2\pi i \ZZ$, induces a map (called by the same name)
$\tr:\CC/(2\pi i \ZZ) \ra\RR$.

With the notation of Equation \eqref{eq:defirKalpha}, note that for
all elements $a,b\in\NN$ that are nonzero, we have $r_{K,\alpha}
(a,b,N) =r_{K}(a)\,r_{K}(b)$ if $N$ is large enough and $r_{K,\alpha}
(a,b,N)=r_{K,\alpha} (b,a,N)$.

Let $\R^+_{N}$ be the restriction to $[0,+\infty[$ of the empirical
measure $\R_{N}$ defined in Equation \eqref{eq:defiRsubN}.  Using
Equation \eqref{eq:defirKalpha} and writing $a=\n(q)$ and $b=\n(p+q)$
in the indices of the first sum below for $p$ and $q$ varying in
$\OOO_K$, since $\tr \circ \log= \ln\circ\n$ and by the
multiplicativity of the norm for the second equality below, and by
Lemma \ref{lem:calculavecp} \eqref{item2:calculavecp} for the third
one, we have
\begin{align*}
\R^+_{N}&=\frac{1}{\psi(N)}\sum_{a,b\in \NN\;:\; 0<a<b\le N^2}\;
r_{K,\alpha}(a,b,N)\,\Delta_{\phi(N)(\ln b-\ln a)}\nonumber
\\&=\frac{1}{\psi(N)}
\sum_{p,q\in \OOO_K\;:\; 0<\n(q)<\n(p+q)\le N^2,\;|p|\leq N^\alpha}\;
\Delta_{\phi(N)\tr\log\frac{p+q}{q}}\nonumber
\\&=\frac{1}{\psi(N)}\sum_{p\in \OOO_K\;:\;0<|p|\leq
N^\alpha}\;\sum_{q\in J_{p,N}}
\Delta_{\phi(N)\tr\log(1+\frac{p}{q})}\,.
\end{align*}
Let us fix for now $p\in \OOO_K$ with $0<|p|\leq N^\alpha$. Let us
consider the measure with finite support
\[
\nu_{p,N}=\sum_{q\in J_{p,N}} \Delta_{\phi(N)\tr\log(1+\frac{p}{q})}\,, 
\]
so that we have
\begin{equation}\label{eq:RNplus}
  \R^+_{N}=\frac{1}{\psi(N)}
  \sum_{p\in \OOO_K:\;0<|p|\leq N^\alpha}\;\nu_{p,N}\,.
\end{equation}
By Equation \eqref{eq:defijpN} and since $j_{p,N}$ is a bijection, we
have
\begin{equation}\label{eq:netoynupN}
  \nu_{p,N}=\sum_{(s,t)\in j_{p,N}^{-1}(J_{p,N})} \Delta_{\phi(N)\tr\log\big(1+
  \frac{p}{N(s\,\frac{p}{|p|}+t\,\frac{v_p}{|v_p|})}\big)}\,.
\end{equation}

Let $f\in C^1_c([0,+\infty[)$ be a $C^1$-smooth function on $[0,+\infty[$
with compact support. As $N\ra+\infty$, for every $(s,t)\in\RR^2$ with
$s^2+t^2\geq \epsilon^2$ (where $\epsilon$ is given by Equation
\eqref{eq:conditiongaeps}),

$\bullet$~ by the expansion of $\log (1+z)$ near $z=0$ since
$\big|s\,\frac{p}{|p|}+ t\,\frac{v_p}{|v_p|}\big| =\sqrt{s^2+t^2}\geq
\epsilon$, so that using Equation \eqref{eq:conditiongaeps} we have
$\Big|\frac{p}{N(s\,\frac{p}{|p|}+ t\, \frac{v_p}{|v_p|})}\Big| \leq
\frac{|p|}{N\epsilon}\leq \frac{1}{N^{1-\alpha-\gamma}}$, which tends
to $0$ as $N\ra+\infty$, and by the linearity of the trace, for the
first equality,
    
$\bullet$~ since $p$ and $v_p$ are orthogonal, for the third equality,

$\bullet$~ by the assumption of Equation \eqref{eq:scalingrange} so
that $\frac{\phi(N)\,|p|^2}{N^2\,\epsilon^2}\leq\frac{\phi(N)}
{N^{2-2\alpha-2\ga}}$ tends to $0$ as $N\ra+\infty$, and by the mean
value inequality, for the fourth equality,

\noindent we have
\begin{align}  
   &f\Big(\phi(N)\tr\log\Big(1+
    \frac{p}{N(s\,\frac{p}{|p|}+t\,\frac{v_p}{|v_p|})}\Big)\Big)
    \nonumber\\=\;& f\Big(\;\frac{\phi(N)}{N}\tr\Big(
    \frac{p}{s\,\frac{p}{|p|}+t\,\frac{v_p}{|v_p|}}\Big)
    +\bigO\Big(\frac{\phi(N)\,|p|^2}{\epsilon^2\,N^2}\Big)\Big)
    \nonumber\\ =\;&f\Big(\;\frac{\phi(N)}{N(s^2+t^2)}\tr\Big(
    p\Big(s\,\frac{\ov{p}}{|p|}+t\,\frac{\ov{v_p}}{|v_p|}\Big)\Big)
    +\bigO\Big(\frac{\phi(N)\,|p|^2}{\epsilon^2\,N^2}\Big)\Big)
    \nonumber\\ =\;&f\Big(\,\frac{2\,\phi(N)\,|p|\,s}{N(s^2+t^2)}+\bigO
    \Big(\frac{\phi(N)\,|p|^2}{\epsilon^2\,N^2}\Big)\Big) \nonumber\\
    =\;&f\Big(\,2\,\frac{\phi(N)}{N}\,|p|\,\frac{s}{s^2+t^2}\,\Big)
    +\bigO\Big(\;\frac{\phi(N)\,|p|^2}{N^2\,\epsilon^2}\,
    \|f'\|_\infty\Big)\,. \label{eq:expressRNplus}
\end{align}

As in Section \ref{sec:unfold}, let us denote by $B(z,r)$ the closed
ball of center $z$ and radius $r$ in $\CC$.  The function 
\[
\wh f=\wh f_{p,N}:(s,t)\mapsto
f\Big(\,2\, \frac{ \phi(N)}{N} \,|p|\,\frac{s}{s^2+t^2}\,\Big)
\]
is well defined and $C^1$-smooth on $B(0,2)\ssm B(0,\frac{\epsilon}
{2})$. An easy computation gives that the supremum norm $\|\,d\wh
f\;\|_\infty$ of its differential on $B(0,2)\ssm B(0,\frac{\epsilon}
{2})$ satisfies
\begin{equation}\label{eq:ccontroldiff}
\|\,d\wh f\;\|_\infty=\bigO\Big(\,\frac{\phi(N)\,|p|}{N\,\epsilon^2}
\,\|f'\|_\infty\Big)\,.
\end{equation}
As $N\ra+\infty$ and uniformly in $p\in\OOO_K$ with $0<|p|\leq
N^\alpha$, using

$\bullet$~ the fact that $\card\{j_{p,N}^{-1}(J_{p,N})\cap
B(0,\epsilon)\}= \bigO(N^2\,\epsilon^2)$ as $N\ra+\infty$, by the same
proof as for Equation \eqref{eq:cardJpN}, for the first equality,

$\bullet$~ Equations \eqref{eq:netoynupN} and
\eqref{eq:expressRNplus}, as well as Equation \eqref{eq:cardJpN} for
dealing with the error term of Equation \eqref{eq:expressRNplus}, for
the second equality,

\noindent
we have

\begin{align}
\nu_{p,N}(f)&=\sum_{(s,t)\in j_{p,N}^{-1}(J_{p,N})\smallsetminus B(0,\epsilon)}
f\Big(\phi(N)\tr\log\Big(1+
\frac{p}{N(s\,\frac{p}{|p|}+t\,\frac{v_p}{|v_p|})}\Big)\Big)+\bigO
\big(N^2\,\epsilon^2\,\|f\|_\infty\big)\nonumber
\\&=\sum_{z\in j_{p,N}^{-1}(J_{p,N})\smallsetminus B(0,\epsilon)}
\wh f\,(z) \quad
+\bigO\big(\phi(N)\,|p|^2\,\epsilon^{-2}\,\|f'\|_\infty\big)+\bigO
\big(N^2\,\epsilon^2\,\|f\|_\infty\big)\,.\label{eq:somsurqun}
\end{align}

Let $\Lambda=\frac{c'_p}{2\,|p|\,N}\,\ZZ+\frac{|v_p|}{N}\,\ZZ \,i$,
which is a $\ZZ$-lattice in the Euclidean space $\RR^2=\CC$, with
fundamental parallelogram $\big[0,\frac{c'_p}{2\,|p|\,N}\big] \times
\big[0,\frac{|v_p|}{N}\big]$. Its diameter $\diam_\Lambda$ and area
$\covol_\Lambda$ satisfy respectively by the left part and the right
part of Equation \eqref{eq:controlcppvp} that
\begin{equation}\label{eq:controldiam}
  \diam_\Lambda=
  \Big(\frac{{c'_p}^2}{4\,|p|^2\,N^2}+\frac{|v_p|^2}{N^2}\Big)^{1/2}
  =\bigO\Big(\frac{|p|}{N}\Big)
\end{equation}
and
\begin{equation}\label{eq:controlcovol}\quad
  \covol_\Lambda=\frac{c'_p\,|v_p|}{2\,|p|\,N^2}
  \in\Big[\frac{1}{2\,c_K\,N^2},\frac{c_K}{2\,N^2}\Big]\,.
\end{equation}
If $N$ is large enough and uniformly in $p\in\OOO_K$ with $0<|p|\leq
N^\alpha$, for all integers $k,\ell\in\ZZ$, since the real number
$t_{p,k}$ given by Lemma \ref{lem:calculavecp}
\eqref{item5:calculavecp} satisfies the inequality $|t_{p,k}| \leq
\frac12$, every point $z= \big(\,\frac{c'_p} {2\,|p|\,N}
\,k,\frac{t_{p,k}\,|v_p|}{N} +\frac{|v_p|}{N} \,\ell \,\big)$ in
$j_{p,N}^{-1}(J_{p,N})\ssm B(0,\epsilon)$ is at distance at most
$\frac{|v_p|}{N}$ from the point $z'= \big(\,\frac{c'_p}{2\,|p|\,N}
\,k,\frac{|v_p|}{N}\,\ell\,\big)$ of $\Lambda$.  By the left part of
Equation \eqref{eq:controlcppvp} and since $\ga < 1-\alpha$, for every
$N$ large enough, we have
\[
\frac{|v_p|}{N}\leq\frac{c_K\,|p|}{N}\leq
\frac{c_K}{N^{1-\alpha}} < \frac{1}{2\,N^\ga}= \frac{\epsilon}{2}\,.
\]
Recalling that we have $j_{p,N}^{-1}(J_{p,N})\subset B\big(-\frac{|p|}
{N},1\big)$, we hence have $z,z'\in B(0,2) \ssm B(0,\frac{\epsilon}
{2})$ for every $N$ large enough. Again by the mean value theorem and
by Equation \eqref{eq:ccontroldiff}, we therefore have
\[
\big|\,\wh f\,(z)-\wh f\,(z')\,\big|=
\bigO\Big(\frac{|p|}{N}\,\|\,d\wh f\;\|_\infty\Big) =
\bigO\Big(\,\frac{\phi(N)\,|p|^2}{N^2\,\epsilon^2}\,\|f'\|_\infty\Big)\,.
\]
Thus by Equations \eqref{eq:somsurqun} and \eqref{eq:cardJpN}, we have
\begin{align*}
\nu_{p,N}(f)&=\sum_{z\in j_{p,N}^{-1}(J_{p,N})\smallsetminus B(0,\epsilon)}
\wh f\,(z')\nonumber \\&\qquad\quad
+\bigO\big(\phi(N)\,|p|^2\,\epsilon^{-2}\,\|f'\|_\infty\big)+\bigO
\big(N^2\,\epsilon^2\,\|f\|_\infty\big)\,.%\label{eq:somsurqdeux}
\end{align*}
The symmetric difference between the set $\Lambda\cap(B^+(0,1)\ssm
B(0,\epsilon))$ and the set of elements $z'$ such that $z\in
j_{p,N}^{-1}(J_{p,N})\ssm B(0,\epsilon)$ is, by the triangle
inequality, contained in the intersection of $\Lambda$ with the
$2\frac{|p|}{N}$-neighbourhood $\N_{2\frac{|p|}{N}}(\partial
B^+(0,1))$ of the boundary of $B^+(0,1)$. By the Gauss counting
argument and by Equation \eqref{eq:controlcovol}, this intersection
has cardinality $\bigO\big(\frac{1}{\covol_\Lambda}
\operatorname{Leb}_\CC\big(\N_{2\frac{|p|}{N}}(\partial B^+(0,1))\big)
\big) =\bigO(N\,|p|)$. Hence
\begin{align}
\nu_{p,N}(f)=&\sum_{z'\in \Lambda\cap(B^+(0,1)\smallsetminus B(0,\epsilon))}
\wh f\,(z') \nonumber\\ & \qquad\quad
+\bigO\big(\phi(N)\,|p|^2\,\epsilon^{-2}\,\|f'\|_\infty\big)+\bigO
\big((N\,|p|+N^2\,\epsilon^2)\,\|f\|_\infty\big)\,.\label{eq:somsurqtrois}
\end{align}
By the usual approximation of two-dimensional integrals by Riemann
sums, we have
\begin{equation}\label{eq:2driem}
\Big|\int_{B^+(0,1)\smallsetminus B(0,\epsilon)} \wh f\;d\Leb_\CC-
\covol_\Lambda\sum_{z'\in \Lambda\cap(B^+(0,1)\smallsetminus B(0,\epsilon))}
\wh f\,(z') \Big|=\bigO\big(\diam_\Lambda\,\|\,d\wh f\;\|_\infty\big)\,.
\end{equation}
By Equations \eqref{eq:controldiam}, \eqref{eq:controlcovol} and
\eqref{eq:ccontroldiff}, we have
\[
\frac{\diam_\Lambda}{\covol_\Lambda}\,\|\,d\wh f\;\|_\infty
=\bigO(\phi(N)|p|^2\epsilon^{-2}\,\|f'\|_\infty)\,.
\]
Therefore Equation \eqref{eq:somsurqtrois} becomes, using Equations
\eqref{eq:controlcovol} and \eqref{eq:2driem} for the first equality
below, and the fact that the area of $B(0,\epsilon)$ is
$\bigO(\epsilon^2)$ and Equation \eqref{eq:controlcppvp} for the
second one,
\begin{align}
  \nu_{p,N}(f)=&\frac{2\,|p|\,N^2}{c'_p\,|v_p|}
  \int_{B^+(0,1)\smallsetminus B(0,\epsilon)} \wh f\;d\Leb_\CC
  \nonumber\\ & \qquad
+\bigO\big(\phi(N)\,|p|^2\,\epsilon^{-2}\,\|f'\|_\infty\big)+\bigO
\big((|p|\,N+N^2\,\epsilon^2)\,\|f\|_\infty\big)\nonumber\\=
&\frac{2\,|p|\,N^2}{c'_p\,|v_p|}
  \int_{B^+(0,1)} \wh f\,d\Leb_\CC
  \nonumber\\ & \qquad
+\bigO\big(\phi(N)\,|p|^2\,\epsilon^{-2}\,\|f'\|_\infty\big)+\bigO
\big((|p|\,N+N^2\,\epsilon^2)\,\|f\|_\infty\big)\,.\label{eq:somsurqquatre}
\end{align}

Let us define $\wt f=\wt f_{p,N}\in C_c([0,+\infty[)$ by $\wt f:u
\mapsto f\big(2\,\frac{\phi(N)}{N}\,|p|\,u\big)$, so that we have $\wh
f(s,t)=\wt f\big(\frac{s}{s^2+t^2}\big)$ for all $(s,t)\in B^+(0,1)
\ssm\{0\}$. Let us now compute the pushforward measure of ${\rm Leb}
_{B^+(0,1)}$ by the map from $\RR^2$ to $\RR$ defined by $(s,t)\mapsto
\frac{s}{s^2+t^2}$. Neglecting sets of measure $0$ and using

$\bullet$~ Equation \eqref{eq:Bplus01} and the symmetry $t\mapsto -t$,
for the first equation,

$\bullet$~ the change of variable (with $s$ fixed)
    $u=\frac{s}{s^2+t^2}\geq 0$ so that $t=\sqrt{\frac{s}{u}-s^2}$ and
    $dt=-\frac{ \sqrt{s} }{ 2\,u^{3/2} \,\sqrt{1-us} }\;du$, for the
    second equation,

$\bullet$~ the fact that $u\in[s,1/s]$ if and only if $s\in[0,\min\{u,
\frac{1}{u}\}]$ and Fubini's theorem, for the third equation,

\noindent
we have
\begin{align}  
  &\int_{B^+(0,1)}\wt f\Big(\,\frac{s}{s^2+t^2}\,\Big)\;ds\,dt  =
  2\int_{s=0}^1\int_{t=0}^{\sqrt{1-s^2}}\wt f\Big(\,\frac{s}{s^2+t^2}
  \,\Big)\;dt\,ds\nonumber\\ =\;& \int_{s=0}^1\int_{u=s}^{1/s} \wt f(u)\;
  \frac{ \sqrt{s} }{u^{3/2} \,\sqrt{1-us} }\;du\,ds=
  \int_{u=0}^{+\infty}\frac{\wt f(u)}{u^{3/2}}
  \int_{s=0}^{\min\{u,\frac{1}{u}\}} \sqrt{\frac{s}{1-us} }\;ds\,du\,.
  \label{eq:pushforw}
\end{align}
Using successively the changes of variable (with $u$ fixed) $\sigma
=us$ and $\theta =\arcsin\sqrt{\sigma}$, and setting $m_u=\min\{1,
u^2\}$, an easy computation gives
\begin{align}
\frac{1}{u^{3/2}}\int_{s=0}^{\min\{u,\frac{1}{u}\}}
\sqrt{\frac{s}{1-us} }\;ds &=\frac{1}{u^3}\int_{\sigma=0}^{m_u}
\sqrt{\frac{\sigma}{1-\sigma} }\;d\sigma =
\frac{1}{u^3}\int_{\theta=0}^{\arcsin(\sqrt{m_u}\,)}
2\sin^2\theta\;d\theta \nonumber\\& = \frac{1}{u^3}
\big(\arcsin(\sqrt{m_u}\,)-\sqrt{m_u}\sqrt{1-m_u}\;\big)\,.
\label{eq:appeararcsin}
\end{align}
Note that $\sqrt{m_u}=\min\{1,u\}$. Consider the function $g:\,]\,0,
  +\infty[\;\ra\RR$ defined by
\begin{equation}\label{eq:defig}
g:u\mapsto \Big\{\begin{array}{ll}
\frac{1}{u^3}(\arcsin(u)-u\sqrt{1-u^2}\,)&\text{if~} u\leq
1\\ \frac{\pi}{2\,u^3}&\text{otherwise.}\end{array}
\end{equation}

\noindent\begin{minipage}{8cm}\begin{center}
\includegraphics[height=4cm]{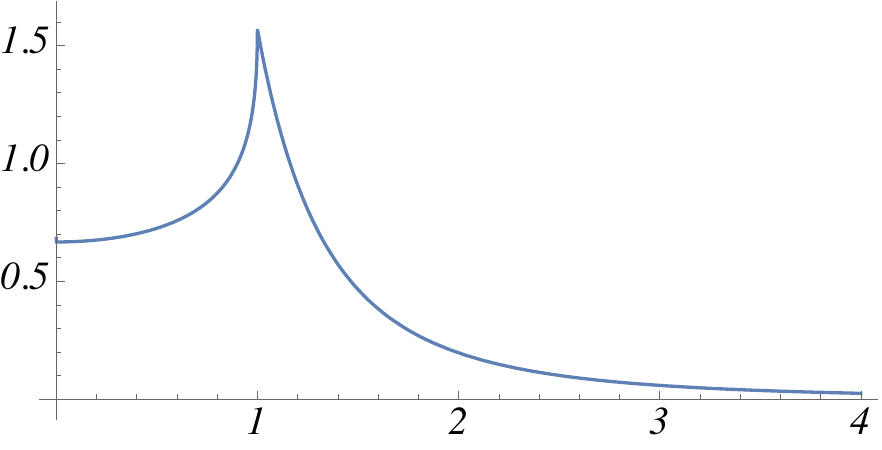}
\end{center}
\end{minipage}
\begin{minipage}{7cm}
  \addtocounter{fig}{1}
  \hbox{\rm Figure \arabic{fig} ~: The graph of the function $g$.}
\end{minipage}

\smallskip
Let us recall the asymptotic expansions near $u=0$ of $\arcsin u=u+
\frac{u^3}{6} +\bigO(u^5)$ and $\sqrt{1-u^2}= 1-\frac{u^2}{2}+
\bigO(u^4)$. Hence the function $g$ extends continuously at $0$ by
$g(0)=\frac{2}{3}$. It is continuous on $[0,+\infty[$, positive, with
upper bound $\|g\|_\infty=g(1)=\frac{\pi}{2}$. It is integrable, and
$C^1$-smooth on $[0,+\infty[$ except at $u=1$, where $g$ is not
differentiable on the left. More precisely, since the derivatives of
both $u\mapsto -u\sqrt{1-u^2}$ and $u\mapsto\arcsin u$ are
Landau-equivalent to $\frac{1}{\sqrt{1-u^2}}$ as $u$ tends to $1$ from
below, we have
\begin{equation}\label{eq:derivgvers1}
  g'(u)=\bigO\Big(\frac{1}{\sqrt{1-u}}\;\Big)\qquad\text{as}\quad u\ra 1^-\,.
\end{equation}
By Equation \eqref{eq:pushforw} (which remains valid) applied with
$\wt f$ the characteristic function $\mathbbm{1}_{[0,1]}$ of $[0,1]$,
we have $\int_0^{1}g(s)\;ds=\int_0^{+\infty}\mathbbm{1}_{[0,1]}(s)\,g(s)
\;ds= 2\int_{s=0}^1\int_{t=\sqrt{s-s^2}}^{\sqrt{1-s^2}}\;dt\,ds=
\frac{\pi}{2}-\frac{\pi}{4}=\frac{\pi}{4}$, hence 
\begin{equation}\label{eq:valintg}
\int_0^{+\infty}g(s)\;ds=\int_0^{1}g(s)\;ds+\int_1^{+\infty}g(s)\;ds=
\frac{\pi}{4}+\frac{\pi}{4}=\frac{\pi}{2}\,.
\end{equation}
Furthermore, we have $g(u)= \frac{2}{3} +\bigO(u^2)$ as $u$ tends to
$0$.  Hence by Equations \eqref{eq:pushforw} and
\eqref{eq:appeararcsin}, and by the definition of $\wt f$, we have
\begin{align*}  
  \int_{B^+(0,1)}\wh f\;d\Leb_\CC&=\int_{B^+(0,1)}\wt f
  \Big(\,\frac{s}{s^2+t^2}\,\Big)\;ds\,dt  =
  \int_{u=0}^{+\infty}\wt f(u)\,g(u)\;du\\&=
  \int_{u=0}^{+\infty} f\Big(2\,\frac{\phi(N)}{N}\,|p|\,u\Big)\,g(u)\;du\,.
\end{align*}
Equation \eqref{eq:somsurqquatre} hence becomes, using the change of
variable $t= 2\,\frac{\phi(N)}{N}\,|p|\,u$ for the second equality,
\begin{align}
  \nu_{p,N}(f)&=\frac{2\,|p|\,N^2}{c'_p\,|v_p|}
  \int_{u=0}^{+\infty} f\Big(2\,\frac{\phi(N)}{N}\,|p|\,u\Big)\,g(u)\;du
  \nonumber\\ & \qquad
+\bigO\big(\phi(N)\,|p|^2\,\epsilon^{-2}\,\|f'\|_\infty\big)+\bigO
(N\,|p|\,\|f\|_\infty)+\bigO(N^2\,\epsilon^2\,\|f\|_\infty)
\nonumber\\
&=\int_{t=0}^{+\infty} f(t) \frac{N^3}{\phi(N)\,c'_p\,|v_p|}
  \,g\Big(\frac{t\,N}{2\,\phi(N)\,|p|}\Big)\;dt
  \nonumber\\ & \qquad
+\bigO\big(\phi(N)\,|p|^2\,\epsilon^{-2}\,\|f'\|_\infty\big)+\bigO
\big(N\,|p|\,\|f\|_\infty)+\bigO(N^2\,\epsilon^2\,\|f\|_\infty\big)\,.
\label{eq:somsurqsix}
\end{align}

In order to prove the main result of Section \ref{sec:linear}, which
is Equation \eqref{eq:caspasplusinfty} below, we will use the
following classical Gauss counting result. For every $\beta'\in[-1,+
\infty[\,$, by for instance the proof of \cite[Lem.~2.10]{Sayous25}
for the first equality and by Equation \eqref{eq:covoldiamsys} for the
second one, as $x\geq 1$ tends to $+\infty$, we have
\begin{align}
  \sum_{p\in\OOO_K:\;0<|p|\leq x} |p|^{\beta'}
  &=\frac{2\,\pi}{\covol_{\OOO_K}}\;\frac{x^{\beta'+2}}{\beta'+2}+
  \bigO_{\beta'}\Big(\frac{1+\diam^2_{\OOO_K}}{\covol_{\OOO_K}}\;
  x^{\beta'+1}\Big) \nonumber
  \\&=\frac{4\,\pi}{\sqrt{|D_K|}}\;\frac{x^{\beta'+2}}{\beta'+2}+
  \bigO_{\beta'}(\sqrt{|D_K|}\,x^{\beta'+1}) =
  \bigO_{\beta',D_K}(x^{\beta'+2})\,.\label{eq:GaussSayous}
\end{align}
Similarly, we have the following analogous result for the case $\beta'
=-2$ that we will only use in Section \ref{sec:mainproof}: As $x\geq
1$ tends to $+\infty$, we have
\begin{align}\label{eq:GaussSayousmoin2}
  \sum_{p\in\OOO_K:\;0<|p|\leq x} \frac{1}{|p|^{2}}=
  \frac{4\,\pi}{\sqrt{|D_K|}}\;\ln x+ \bigO_{D_K}(1)\,.
\end{align}

For every $N\in\NN\ssm\{0\}$, recalling the definition of the function
$g$ in Equation \eqref{eq:defig}, let us define a function $\Theta_N:
[0,+\infty[\;\ra[0,+\infty[$ by
\begin{equation}\label{eq:defiThetaN}
\Theta_N:t\mapsto \frac{N^3}{\psi(N)\,\phi(N)}
\sum_{p\in \OOO_K:\;0<|p|\leq N^\alpha} \frac{1}{c'_p\,|v_p|}
\,g\Big(\frac{t\,N}{2\,\phi(N)\,|p|}\Big) \,.
\end{equation}
By the assumption  $\ga < \frac{1-\alpha}{2}$  in
Equation \eqref{eq:conditiongaeps}, we have $1+3\alpha\leq
2+2\alpha- 2\ga$. By Equations \eqref{eq:RNplus} and
\eqref{eq:somsurqsix}, by using Equation \eqref{eq:GaussSayous} with
$\beta'=2,1,0$ and $x= N^\alpha$ in order to control the three error
terms in Equation \eqref{eq:somsurqsix}, since $\epsilon = N^{-\ga}$
and since $1+3\alpha\leq 2+2\alpha-2\ga$, we hence have
\begin{align}
  \R^+_{N}(f)&=\frac{1}{\psi(N)}
  \sum_{p\in \OOO_K:\;0<|p|\leq N^\alpha}\;\nu_{p,N}(f)=
  \int_{t=0}^{+\infty} f(t)\,\Theta_N(t)\;dt \nonumber\\
  &\qquad+\bigO\Big(\frac{\phi(N)\,N^{4\alpha+2\ga}}{\psi(N)}\,
  \|f'\|_\infty\Big)+\bigO\Big(\frac{N^{2+2\alpha-2\ga}}{\psi(N)}\,
  \|f\|_\infty\Big)\,.
\label{eq:caspasplusinfty}
\end{align}

\section{The main result and its proof}
\label{sec:mainproof}

Before stating our main result Theorem \ref{theo:sumsquaremain}, let
us give the mild restrictions on the scaling function $\phi$ that we
will use. We keep the notation $\alpha\in\;]0,\frac12[$ of the
introduction. We assume in this section that the limits
\[
\lambda_\phi =
\lim_{N\ra+\infty}\;\frac{\phi(N)}{N^{1-\alpha}},\;
\lambda'_\phi =\lim_{N\ra+\infty}\;
\frac{\phi(N)}{N^{1-\frac{\alpha}{2}}},\;
\lambda''_\phi= \lim_{N\ra+\infty}\;\frac{\phi(N)}{N}
\;\text{and}\;
\lambda'''_\phi= \lim_{N\ra+\infty}\;\frac{\phi(N)}{N^{1+\frac{\alpha}{2}}}
\]
exist in $[0,+\infty]$, and that there exists $\ga>0$ such that
Equations \eqref{eq:conditiongaeps} and \eqref{eq:scalingrange} still
hold. When $\phi$ is a power scaling $N\mapsto N^\beta$ as in the
introduction and if $\beta\in\;]0,2-2\alpha[\,$, these assumptions are
satisfied if and only if $\ga\in\;]0,\min\{\frac{1-2\alpha}{2},
1-\alpha-\frac{\beta}{2}\}[\,$.

For all $k\in\NN$ and $N\in\NN\ssm\{0\}$, let us define
\begin{equation}\label{eq:defiSN}
  S_{N,k}=\frac{|D_K|\,(k+1)}{4\,\pi}\sum_{p\in \OOO_K:\;0<|p|\leq N^\alpha}
  \frac{|p|^k}{c'_p\,|v_p|}\,.
\end{equation}
By the right part of Equation \eqref{eq:controlcppvp} and by Equation
\eqref{eq:GaussSayous} with $\beta'=k-1$ and $x=N^\alpha$, we have
\begin{align*}
S_{N,k}&\leq\frac{c_K\;|D_K|\,(k+1)}{4\,\pi}
\sum_{p\in \OOO_K:\;0<|p|\leq N^\alpha}|p|^{k-1}\\&=
c_K\;\sqrt{|D_K|}\;N^{(k+1)\alpha}+
\bigO_k(c_K\,|D_K|^{\frac{3}{2}}\,N^{k\,\alpha})\,.
\end{align*}
With the similarly obtained lower bound, there hence exists a constant
$c'_{K,k}>0$ such that for every $N\in\NN\ssm\{0\}$, we have
\begin{equation}\label{eq:controlSN}
  \frac{1}{c'_{K,k}}\;N^{(k+1)\alpha} \leq S_{N,k}\leq
  c'_{K,k}\;N^{(k+1)\alpha}\,.
\end{equation}
Note that when $D_K\equiv 0\!\!\mod 4$, we have more precisely by
Equations \eqref{eq:valsppvpdiscnul} and \eqref{eq:GaussSayous} with
$\beta'=k-1$ and $x=N^\alpha$ that
\begin{equation}\label{eq:defiSNcasdiscnul}
  S_{N,k}=\frac{\sqrt{|D_K|}\,(k+1)}{4\,\pi}
  \sum_{p\in \OOO_K:\;0<|p|\leq N^\alpha}|p|^{k-1}
  =N^{(k+1)\alpha}+\bigO_k(|D_K|N^{k\,\alpha})\,.
\end{equation}

We now define the measures $m_{\phi}$ that will appear as asymptotic
pair correlation measures in Theorem \ref{theo:sumsquaremain}. If
$\lambda_\phi\in\;]0,+\infty[\,$, we define an even function
$\rho_{1-\alpha}:\RR\ra[0,+\infty[$ by
\begin{equation}\label{eq:defpaircorrfunctgene}
\rho_{1-\alpha}(t)=\left\{\begin{array}{ll}
\frac{8\,\pi}{3\,|D_K|} &\!\!\text{if } t=0\smallskip\\
\frac{8\,\pi\,\lambda_\phi^3}{|D_K|\,t^3}
\Big(\arcsin\big(\frac{t}{2\lambda_\phi}\big)-
\frac{t}{2\lambda_\phi}\big(1-\frac{t^2}{2\lambda_\phi^2}\big)
\big(1-\frac{t^2}{4\lambda_\phi^2}\big)^{\frac{1}{2}}\,\Big)
&\!\!\text{if } 0<t\leq 2\lambda_\phi \\
\frac{4\,\pi^2\,\lambda_\phi^3}{|D_K|\,t^3}
&\!\!\text{if } t> 2\lambda_\phi\,.
\end{array}\right.
\end{equation}
It is easy to see that $\rho_{1-\alpha}$ is bounded, continuous,
positive and piecewise real analytic. See Equation
\eqref{eq:defpaircorrfunct} and the graph of $\rho_{1-\alpha}$ when
$\lambda_\phi=1$ in Figure 1 of the Introduction (corresponding to the
power scaling $\phi:N\mapsto N^{1-\alpha}$). Let
\[
m_{\phi}=
\begin{cases}
\frac{4\,\pi^2}{|D_K|}\;\Delta_0
&\textrm{if } \lambda_\phi=0,\\
\rho_{1-\alpha}\,\Leb_\RR
  &\textrm{if } \lambda_\phi\in\;]0,+\infty[\,,\\
\frac{8\,\pi}{3\,|D_K|}\;\Leb_\RR
&\textrm{if } \lambda_\phi=+\infty\,.
\end{cases}
\]

\btheo\label{theo:sumsquaremain}
For every $A\geq 1$ and for every $f\in C^1_c(\RR)$  with support contained
in $[-A,A]$, as $N\ra+\infty$, we have
\begin{align*}
&\R_N(f)=\int_\RR f(t)\,dm_\phi(t)+\\&\begin{cases}
    \!\!\!\begin{array}[t]{l}
      \;\;\;\bigO\big(\big(\frac{\phi(N)}{N^{1-\alpha}}\big)^{\frac{1}{2}}\,
      \|f'\|_\infty\big)\\+
      \bigO\big(\max\big\{\frac{\phi(N)}{N^{1-\alpha}},\;
  \frac{1}{N^{2\ga}}\big\}\,\|f\|_\infty\big)\end{array}
    &\!\!\!\begin{array}[t]{l}\textrm{if } \lambda_\phi=0\\
    \text{and } \psi(N)=N^2\,S_{N,1},\end{array}\medskip\\
\!\!\!\begin{array}[t]{l}
      \;\;\;\bigO\big(N^{-\alpha}\,\|f'\|_\infty\big)+
      \bigO\big(\max\big\{\frac{A}{N^{\frac{\alpha}{4}}},
  \frac{1}{N^{2\ga}}\big\}\,\|f\|_\infty\big)\\+\bigO\big(
      A\big|\,\frac{N^{1-\alpha}}{\phi(N)}-\frac{1}{\lambda_\phi}\,
      \big|^{\frac{1}{2}}\big|\ln\big|\,\frac{N^{1-\alpha}}{\phi(N)}-
  \frac{1}{\lambda_\phi}\,\big|\,\big|\,\|f\|_\infty\big)\end{array}
    &\!\!\!\begin{array}[t]{l}\textrm{if } \lambda_\phi\in\;]0,+\infty[\,,
    \;D_K\equiv 0\!\!\!\!\mod 4,\\ \text{and }
    \psi(N)=\frac{N^{3+\alpha}}{\phi(N)},\end{array}
    \medskip\\\!\!\!\begin{array}[t]{l}
    \;\;\;\bigO\big(N^{\frac{5\alpha-2}{4}}\,\|f'\|_\infty\big)\\+
  \bigO\big(\max\big\{A^3\big(\frac{N^{1-\alpha}}{\phi(N)}\big)^\frac{2}{3}
  ,\;N^{\frac{5\alpha-2}{4}}\big\}\,\|f\|_\infty\big)\end{array}
    &\!\!\!\begin{array}[t]{l}\textrm{if } \lambda_\phi=+\infty,
    \lambda'_\phi<+\infty, \alpha<\frac{2}{5}\\\text{and }
    \psi(N)=\frac{N^3\,S_{N,0}}{\phi(N)},\end{array}\medskip\\
\bigO\big(\frac{\phi(N)}{N^{1+\frac{\alpha}{2}}}\|f'\|_\infty+
\frac{N^{1-\frac{\alpha}{2}}}{\phi(N)}\|f\|_\infty\big)
    &\!\!\!\begin{array}[t]{l}\textrm{if } \lambda'_\phi=+\infty,\;
\lambda''_\phi=0,\;\alpha\leq\frac{2}{11},\\\text{and }
\psi(N)=\frac{N^3\,S_{N,0}}{\phi(N)},\end{array}\medskip\\
\bigO\big(\frac{\phi(N)}{N^{1+\frac{\alpha}{2}}}
  \,(\|f'\|_\infty+A^3\,\|f\|_\infty)\big)
  &\!\!\!\begin{array}[t]{l}\textrm{if } \lambda''_\phi>0,\;
  \lambda'''_\phi=0,\; \alpha\leq\frac{2}{11},\\
     \text{and } \psi(N)=\frac{N^3\,S_{N,0}}{\phi(N)}\,.\end{array}
\end{cases}
\end{align*}

\etheo

Note that we have $\lambda'''_\phi\leq\lambda''_\phi\leq\lambda'_\phi
\leq\lambda_\phi$ for the extended order on $[0,+\infty]$. Hence if
$\lambda''_\phi=+\infty$, then $\lambda'_\phi= \lambda_\phi =+\infty$
and if $\lambda_\phi<+\infty$, then $\lambda'_\phi = \lambda''_\phi
=\lambda'''_\phi =0$. Hence, for instance when $\alpha\leq
\frac{1}{6}$, the list of cases in the above Theorem
\ref{theo:sumsquaremain} is complete, except that the case
$\lambda_\phi\in\;]0,+\infty[\,, \;D_K\notequiv 0\!\!\mod 4$ and the
case $\lambda'''_\phi>0$ are missing.  The first one should be handled
similarly, though the computational complexity seems to be much
higher. For the second one, we refer to Section \ref{rem:casBetageq2}.

\medskip
\dem The pushforward of a measure $\mu$ by a mapping $h$ is denoted by
$h_*\mu$. We denote by $\sg:\RR\ra\RR$ the change of sign map
$t\mapsto -t$. By the change of variables $(a,b)\mapsto (b,a)$ in the
summation of Equation \eqref{eq:defiRsubN}, we have $\sg_*\R_{N}=
\R_{N}$ for every $N\in\NN\ssm\{0\}$. Since the above measures
$m_\phi$ are invariant under $\sg$, we hence only have to prove
Theorem \ref{theo:sumsquaremain} where the empirical measure $\R_{N}$
is replaced by its restriction $\R^+_{N}$ to $[0,+\infty[\,$.
    
We fix throughout this proof $A\geq 1$ and $f\in C^1_c([0,+\infty[)$ a
$C^1$-smooth function on $[0,+\infty[$ with compact support contained
in $[0,A]$. We consider $N\in\NN\ssm\{0\}$ large enough.  

\medskip
We now separate the proof of Theorem \ref{theo:sumsquaremain} into the
five cases appearing in its statement, corresponding to an increasing
scaling. For reasons that will become clear, we will prove Case
\hyperlink{cas4}{4} after Case \hyperlink{cas5}{5}.

\medskip
\noindent{\bf Case \hypertarget{cas1}{1}. } Let us assume that
$\lambda_\phi={\displaystyle\lim_{N\ra+\infty}}
\frac{\phi(N)}{N^{1-\alpha}} =0$.

Let us take $\psi$ to be the function $N\mapsto N^2S_{N,1}$, with $S_{N,1}$
defined in Equation \eqref{eq:defiSN} for $k=1$. Let $\ga$ be any
element of $]0,\frac{1-2\alpha}{2}[\,$. The assumptions
\eqref{eq:conditiongaeps} and \eqref{eq:scalingrange} are satisfied
since $\ga\leq\frac{1-\alpha}{2}$ and ${\displaystyle\lim_{N\ra+\infty}}
\frac{\phi(N)}{N^{1-\alpha}} =0$. Hence we may apply the results of
Section \ref{sec:linear}.

We are going to prove that in Case \hyperlink{cas1}{1}, as $N\ra
+\infty$, the measure $\Theta_N\Leb_{[0,+\infty[}$ on $[0,+\infty[$,
with $\Theta_N$ defined in Equation \eqref{eq:defiThetaN},
weak-star converges to the Dirac mass $\frac{2\,\pi^2}{|D_K|}\,
\Delta_0$ at $0$ with weight $\frac{2\,\pi^2}{|D_K|}$. Below are the
graphs of $\Theta_N$ for various $N$, for the power scaling
$\phi(N)=N^\beta$, in the Gaussian case $K=\QQ[i]$, with $\alpha=0.15$
and $\beta=0.8<1-\alpha$.
\begin{center}
    \includegraphics[width=14cm]{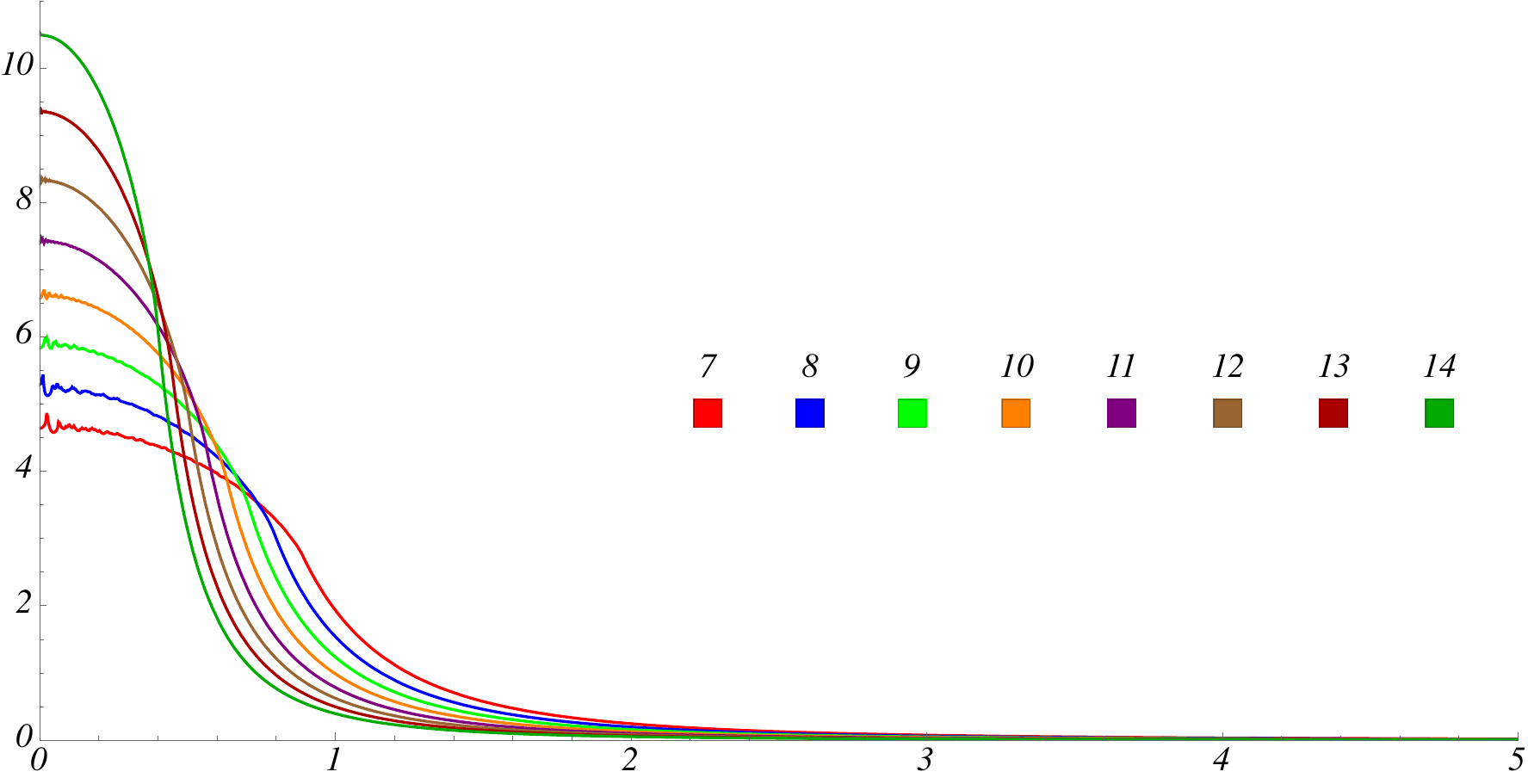}
\end{center}\vspace*{-0.5cm}
\begin{center}
\addtocounter{fig}{1} {\small Figure \arabic{fig}~: Graph of
  $\Theta_N$ for $N=10^m$ with $m\in\{7,8,9,10,11,12,13,14\}$.}
\end{center}
Recall that by Equation \eqref{eq:valintg}, the positive function $g$
is integrable over $[0,+\infty[\,$ with $\int_0^{+\infty}g(s)\;ds=
\frac{\pi}{2}$.  For every $N\in\NN\ssm\{0\}$, by Equation
\eqref{eq:defiThetaN} for the first equality, by the change of
variable $s=\frac{t\,N} {2\,\phi(N)\,|p|}$ for the second equality,
since by Equation \eqref{eq:defiSN} for $k=1$ we have
\[
\sum_{p\in \OOO_K:\;0<|p|\leq N^\alpha}
\frac{|p|}{c'_p\,|v_p|}=\frac{4\,\pi\,S_{N,1}}{2\,|D_K|}\,,
\]
and by the definition of $\psi(N)$ for the last equality, we have
\begin{align}
\int_0^{+\infty}\Theta_N(t)\;dt&=\frac{N^3}{\psi(N)\,\phi(N)}
\sum_{p\in \OOO_K:\;0<|p|\leq N^\alpha} \frac{1}{c'_p\,|v_p|}
\int_0^{+\infty}g\Big(\frac{t\,N}{2\,\phi(N)\,|p|}\Big)\;dt
\nonumber\\&=\frac{2\,N^2}{\psi(N)} \sum_{p\in \OOO_K:\;0<|p|\leq
  N^\alpha} \frac{|p|}{c'_p\,|v_p|} \int_0^{+\infty}g(s)\;ds =
\frac{2\,\pi^2}{|D_K|}\,.\label{eq:valintThetaN}
\end{align}

Let us now prove that the function $\Theta_N$ converges uniformly on
compact subsets of $]0,+\infty[$ to $0$ as $N\ra+\infty$. With the
above centered equation, this will prove that the measure $\Theta_N
\Leb_{[0, +\infty[}$ weak-star converges to the Dirac mass
$\frac{2\,\pi^2}{|D_K|}\,\Delta_0$ as $N\ra+\infty$.

For every $N\in\NN\ssm\{0\}$, let us define
\[
\eta_N=\Big(\frac{\phi(N)}{N^{1-\alpha}}\Big)^{\frac{1}{2}}\,,
\]
noting that under the assumption of Case \hyperlink{cas1}{1}, we have
$\eta_N\ra0$ as $N\ra +\infty$.

Let $t\in[\eta_N,+\infty[$. For every $p\in\OOO_K$ such that $0<|p|
\leq N^\alpha$, we have in particular $\frac{t\,N}{2\,\phi(N)\,|p|}
\geq \frac{\eta_N}{2} \frac{N^{1-\alpha}}{\phi(N)}=\frac{1}{2\,\eta_N}$,
which tends to $+\infty$ as $N\ra+\infty$. In particular, if $N$ is
large enough, we have $u=\frac{t\,N} {2\,\phi(N)\,|p|} \geq 1$. Thus
if $N$ is large enough, respectively

$\bullet$~ by Equation \eqref{eq:defiThetaN} and Equation
\eqref{eq:defig} when $u\geq1$,

$\bullet$~ since $\psi(N)=N^2S_{N,1}$ and by the definition of
$S_{N,3}$ in Equation \eqref{eq:defiSN} for $k=3 $,

$\bullet$~ by Equation \eqref{eq:controlSN} with $k=1$ and $k=3$,

\noindent we have
\begin{align}
  \Theta_N(t)&=\frac{N^3}{\psi(N)\,\phi(N)}
  \sum_{p\in \OOO_K:\;0<|p|\leq N^\alpha} \frac{\pi}{2\,c'_p\,|v_p|} \;
  \frac{8\,\phi(N)^3\,|p|^3}{t^3\,N^3}\nonumber\\&=
  \frac{N^3\,\pi}{N^2\,S_{N,1}\,\phi(N)}\;\frac{4\,\pi\,S_{N,3}}{|D_K|\,4}\;
  \frac{8\,\phi(N)^3}{2\,t^3\,N^3}=
  \frac{4\,\pi^2\,\phi(N)^2\,S_{N,3}}{|D_K|\,t^3\,N^2\,S_{N,1}}\nonumber
  \\&=\frac{1}{t^3}\bigO\Big(\frac{\phi(N)^2\,N^{4\alpha}}
    {N^2\,N^{2\alpha}}\Big)=\frac{1}{t^3}\bigO\Big(
    \Big(\frac{\phi(N)}{N^{1-\alpha}}\Big)^2\Big)\,.
    \label{eq:controlThetaNcas1}
\end{align}
Hence under the assumption of Case \hyperlink{cas1}{1}, the function
$\Theta_N$ converges uniformly on compact subsets of $]0,+\infty[$ to
$0$ as $N\ra+\infty$. But we will need more information on the error
terms.
    
As $N\ra+\infty$, respectively by the additivity of the integral, by
the mean value theorem, by Equation \eqref{eq:valintThetaN} (that gives
$\int_0^{+\infty}\Theta_N(t)\;dt= \frac{2\,\pi^2}{|D_K|}$), by Equation
\eqref{eq:controlThetaNcas1}, and by the value $\eta_N=\Big(
\frac{\phi(N)}{N^{1-\alpha}} \Big)^{\frac{1}{2}}$, we have
\begin{align}
  &\int_0^{+\infty} \Theta_N(t)\,f(t)\;dt=
  \int_0^{\eta_N}\Theta_N(t)\,f(t)\;dt+
  \int_{\eta_N}^{+\infty}\Theta_N(t)\,f(t)\;dt\nonumber\\=\;&
  \int_0^{\eta_N}\Theta_N(t)\,\big(f(0)+
  \bigO(\eta_N\,\|f'\|_\infty)\big)\;dt+
  \int_{\eta_N}^{+\infty}\Theta_N(t)\,f(t)\;dt\nonumber\\=\;&
  \frac{2\,\pi^2}{|D_K|} \,f(0)+\bigO(\eta_N\,\|f'\|_\infty)+
  \int_{\eta_N}^{+\infty}\Theta_N(t)\,(f(t)-f(0))\;dt\nonumber\\=\;&
  \frac{2\,\pi^2}{|D_K|} \,f(0)+\bigO(\eta_N\,\|f'\|_\infty)+
  \int_{\eta_N}^{+\infty}\frac{dt}{t^3}\;\bigO\Big(
  \Big(\frac{\phi(N)}{N^{1-\alpha}}\Big)^2
  \,\|f\|_\infty\Big)\nonumber\\=\;&\frac{2\,\pi^2}{|D_K|} \,f(0)+
    \bigO\Big(\Big(\frac{\phi(N)}{N^{1-\alpha}}\Big)^{\frac{1}{2}}\,
    \|f'\|_\infty\Big)+\bigO\Big(\,\frac{\phi(N)}{N^{1-\alpha}}
    \;\|f\|_\infty\Big)\,.\label{eq:controlThetaNfcas1}
\end{align}

As $N\ra+\infty$, since $\psi(N)=N^2\,S_{N,1}$, by Equation
\eqref{eq:controlSN} with $k=1$, and since we have $2-2\alpha-2\ga\geq
1-\alpha$ as $\ga\leq \frac{1-\alpha}{2}$ by Equation
\eqref{eq:conditiongaeps}, we have
\begin{equation}\label{eq:controlerror1cas1}
  \frac{\phi(N)\,N^{4\alpha+2\ga}}{\psi(N)}=
  \frac{\phi(N)\,N^{4\alpha+2\ga}}{N^2\,S_{N,1}}=
  \bigO\Big(\frac{\phi(N)\,N^{4\alpha+2\ga}}{N^2\,N^{2\alpha}}\Big)=
  \bigO\Big(\frac{\phi(N)}{N^{1-\alpha}}\Big)=
  \bigO\Big(\Big(\frac{\phi(N)}{N^{1-\alpha}}\Big)^{\frac{1}{2}}\Big)\,.
\end{equation}
Similarly, we have
\begin{equation}\label{eq:controlerror2cas1}
  \frac{N^{2+2\alpha-2\ga}}{\psi(N)}= \frac{N^{2+2\alpha-2\ga}}{N^2\,S_{N,1}}=
  \bigO\Big(\frac{N^{2+2\alpha-2\ga}}{N^{2+2\alpha}}\Big)=
  \bigO\Big(\frac{1}{N^{2\ga}}\Big)\,.
\end{equation}
Therefore, using in Equation \eqref{eq:caspasplusinfty} the three Equations
\eqref{eq:controlThetaNfcas1}, \eqref{eq:controlerror1cas1} and
\eqref{eq:controlerror2cas1}, we have
\begin{align}
  \R^+_{N}(f)&=\frac{2\,\pi^2}{|D_K|}\,f(0)+
  \bigO\Big(\Big(\frac{\phi(N)}{N^{1-\alpha}}\Big)^{\frac{1}{2}}\,
  \|f'\|_\infty\Big)+\bigO\Big(\max\Big\{\frac{\phi(N)}{N^{1-\alpha}}
  ,\;\frac{1}{N^{2\ga}}\Big\}\,\|f\|_\infty\Big)\,.
\label{eq:cas1plus}
\end{align}
By symmetry, the restriction of $\R_N$ to $]-\infty,0]$ also
contributes $\frac{2\,\pi^2}{|D_K|}\Delta_0$ to the limit measure,
with an error term as in Equation \eqref{eq:cas1plus} for every $f\in
C_c(\RR)$.  This proves the first case of Theorem
\ref{theo:sumsquaremain}.

When $\phi:N\ra N^\beta$ is a power function, the assumption of Case
\hyperlink{cas1}{1} on $\phi$ are satisfied if and only if
$\beta\in\;]0,1-\alpha[\,$. If furthermore $D_K\equiv 0\!\!\mod 4$,
then by Equation \eqref{eq:defiSNcasdiscnul} with $k=1$, as
$N\ra+\infty$, we have $\psi(N)=N^2\,S_{N,1}\sim N^{2+2\alpha}$. This
proves the first case of Theorem \ref{theo:mainintro}.

\medskip
\noindent{\bf Case \hypertarget{cas2}{2}. } Let us assume that
$\lambda_\phi= {\displaystyle\lim_{N\ra+\infty}} \frac{\phi(N)}
{N^{1-\alpha}} \in\;]0,+\infty[\,$ and that $D_K\equiv 0\!\!\mod 4$.
    
Let us take $\psi$ to be the function $N\mapsto \frac{N^{3+\alpha}}
{\phi(N)}$. Let $\ga$ be any element of $]0,\frac{1-2\alpha}{2}[\,$.
The assumptions \eqref{eq:conditiongaeps} and \eqref{eq:scalingrange}
are satisfied since $\ga<\frac{1-\alpha}{2}$ and
${\displaystyle\lim_{N\ra+\infty}}\frac{\phi(N)}{N^{1-\alpha}}<+\infty$.
Hence we may apply the results of Section \ref{sec:linear}.
  
Below are the graphs of $\Theta_N$ for various $N$, for the power
scaling $\phi:N\mapsto N^\beta$, in the Gaussian case $K=\QQ[i]$, with
$\alpha=0.15$ and $\beta =1-\alpha=0.85$.
\begin{center}
    \begin{overpic}[width=14cm]{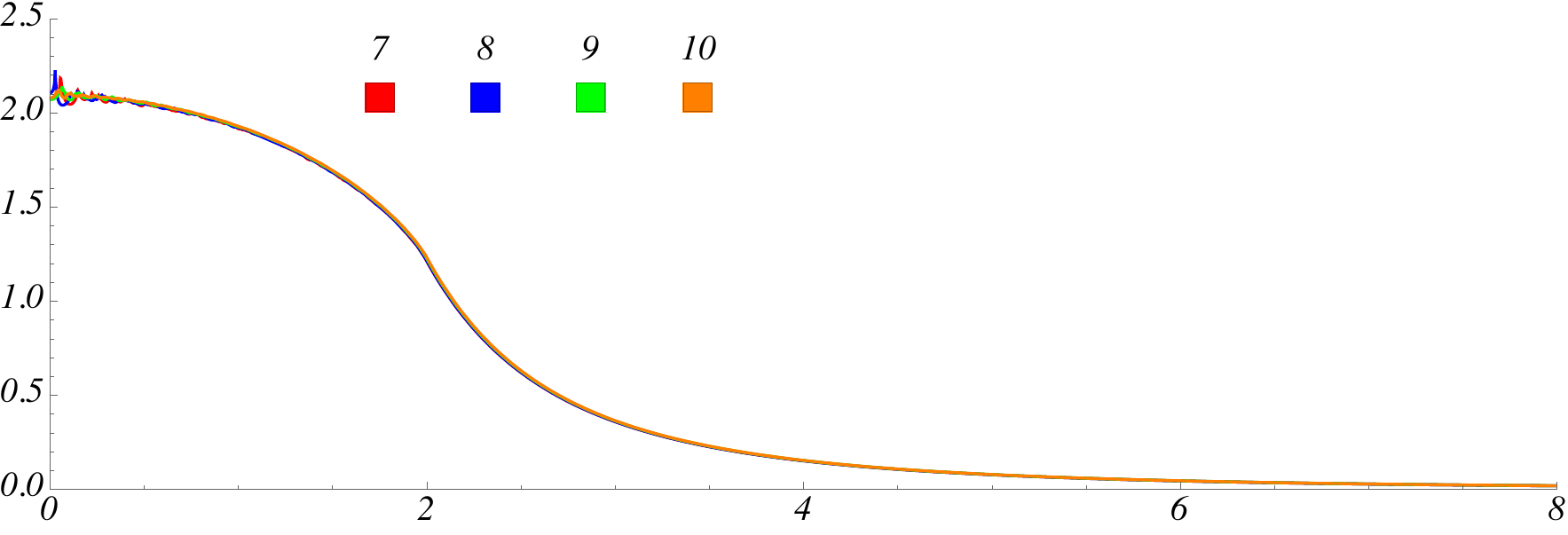}
 \put(50,10){\color{gray}
\frame{\includegraphics[scale=.35]{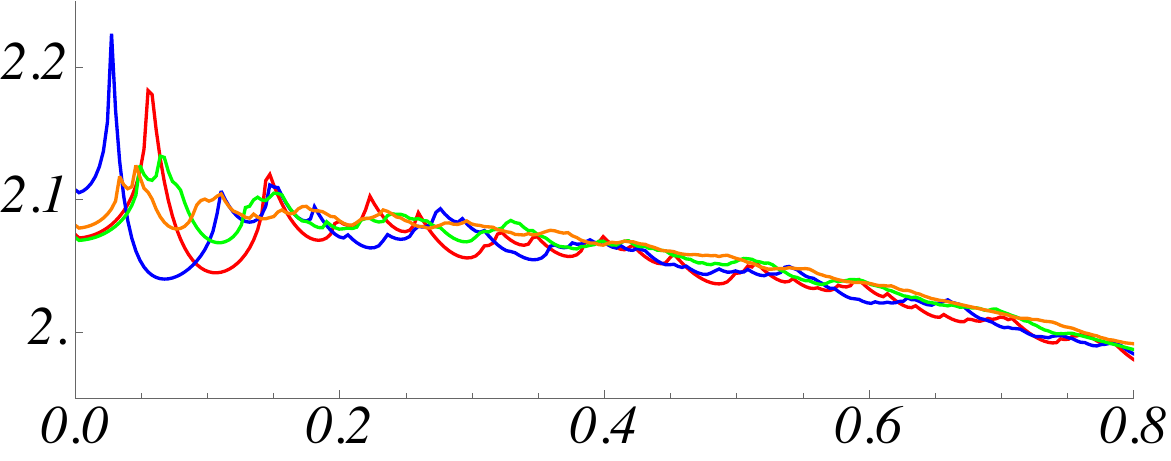}}}
    \end{overpic}
  \end{center}\vspace*{-0.3cm} \begin{center}
\addtocounter{fig}{1} {\small Figure \arabic{fig}~: Graph of
  $\Theta_N$ for $N=10^m$ where $m$ is $7$ (red), $8$ (blue), $9$ (green),
  $10$ (orange).}
\end{center}

Let $t\in[0,A]$ and $N\in\NN\ssm\{0\}$. Since $D_K\equiv 0\!\!\mod 4$
under the assumptions of Case \hyperlink{cas2}{2}, since for every
$p\in\OOO_K\ssm\{0\}$ we then have $c'_p\,|v_p|=\sqrt{|D_K|}\;|p|$ by
Equation \eqref{eq:valsppvpdiscnul}, and since $\psi(N)=
\frac{N^{3+\alpha}}{\phi(N)}$, Equation \eqref{eq:defiThetaN} becomes
\begin{align}
\Theta_N(t)&=\frac{N^3}{\sqrt{|D_K|}\;\psi(N)\,\phi(N)}
\sum_{p\in \OOO_K:\;0<|p|\leq N^\alpha} \frac{1}{|p|}
\,g\Big(\frac{t\,N}{2\,\phi(N)\,|p|}\Big)\nonumber\\&=
\frac{1}{\sqrt{|D_K|}\,N^{2\alpha}}
\sum_{p\in \OOO_K:\;0<|p|\leq N^\alpha} \frac{N^\alpha}{|p|}
\,g\Big(\frac{N^{1-\alpha}}{\phi(N)}\,\frac{t\,N^\alpha}{2\,|p|}\Big)\,.
\label{eq:nettoyageThetaNcas2}
\end{align}
%Note that since $g(0)=\frac{2}{3}$ and by Equation
%\eqref{eq:GaussSayous} with $\beta'=-1$, as $N\ra+\infty$, we have
%\begin{equation}\label{eq:valTheta_Ndezero}
%\Theta_N(0)=\frac{2}{3\,\sqrt{|D_K|}\,N^{\alpha}}
%\sum_{p\in \OOO_K:\;0<|p|\leq N^\alpha} \frac{1}{|p|}
%=\frac{8\pi}{3\,|D_K|}+
%\bigO\big(\frac{1}{N^{\alpha}}\big)\,.
%\end{equation}

Let us define a decomposition of the summation in Equation
\eqref{eq:nettoyageThetaNcas2} corresponding to the subdivision
$[0,+\infty[ \;=[0,1[\;\cup\,[1,+\infty[$ where the function $g$ given
by Equation \eqref{eq:defig} has two different expressions, plus some
safety zone around $1$. Let us define
\begin{equation}\label{eq:defietaetap}
  \varepsilon_N=\max\Big\{8{\lambda_\phi}^2\,
  \Big|\,\frac{N^{1-\alpha}}{\phi(N)}-\frac{1}{\lambda_\phi}\,\Big|,\;
  \frac{1}{N^{\frac{\alpha}{4}}}\Big\}>0\,.
\end{equation}
Note that under the assumptions of Case \hyperlink{cas2}{2}, we have
${\displaystyle \lim_{N\ra+\infty}}\,\varepsilon_N=0$. With the usual
convention on empty sums, let
\begin{align}
&\Theta_N^-(t)=\frac{1}{\sqrt{|D_K|}\;N^\alpha}
  \sum_{p\in \OOO_K:\;0<|p|\leq N^\alpha,\;|p|<\frac{t}{2\lambda_\phi}\,N^\alpha}
  \frac{1}{|p|}\,g\Big(\frac{t\,N}{2\,\phi(N)\,|p|}\Big)\,,
  \label{eq:defThetaNmoin}\\&\label{eq:defThetaNzero}
\Theta_N^0(t)=\frac{1}{\sqrt{|D_K|}\;N^\alpha}
\sum_{\substack{p\in \OOO_K:\;0<|p|\leq N^\alpha\\
    \frac{t}{2\lambda_\phi}\,N^\alpha\leq
  |p|<\frac{t+\varepsilon_N}{2\lambda_\phi}\,N^\alpha}}
  \frac{1}{|p|}\,g\Big(\frac{t\,N}{2\,\phi(N)\,|p|}\Big)\,,
  \\\text{and}\quad&
\Theta_N^+(t)=\frac{1}{\sqrt{|D_K|}\,N^{\alpha}}
\sum_{p\in \OOO_K:\;\frac{t+\varepsilon_N}{2\lambda_\phi}\,N^\alpha\leq|p|\leq N^\alpha}
\frac{1}{|p|}
\,g\Big(\frac{N^{1-\alpha}}{\phi(N)}\,\frac{t\,N^\alpha}{2\,|p|}\Big)\,.
\label{eq:defThetaNplu}
\end{align}
So that by Equation \eqref{eq:nettoyageThetaNcas2} we have
\begin{equation}\label{eq:ThetaNsumpluzeromoi}
\Theta_N(t)=\Theta_N^-(t)+\Theta_N^0(t)+\Theta_N^+(t)\,.
\end{equation}

\medskip\noindent {\bf Step \hypertarget{step1}{1} of Case
  \hyperlink{cas2}{2}. } Let us first estimate $\Theta_N^-(t)$ as
$N\ra+\infty$ uniformly in $t\in[0,+\infty[\,$.

\medskip
Since $\frac{t\,N^{\alpha}}{2\,\lambda_\phi\,|p|}> 1$ (and in
particular $t>0$) whenever the index $p\in\OOO_K$ occurs in the
summation defining $\Theta_N^-(t)$, for $N$ large enough, under the
assumptions of Case \hyperlink{cas2}{2}, we also have
$\frac{t\,N}{2\,\phi(N)\,|p|} =\frac{N^{1-\alpha}}{\phi(N)}\,
\frac{t\,N^\alpha}{2\,|p|}> 1$. By the value $g(u)=\frac{\pi}{2\;u^3}$
on $u\in[1,+\infty[$ given by Equation \eqref{eq:defig}, and by
Equation \eqref{eq:GaussSayous} applied with $\beta'=2$ and
$x=\min\big\{N^\alpha,\frac{t}{2\lambda_\phi}N^\alpha\big\}$, Equation
\eqref{eq:defThetaNmoin} gives, as $N\ra+\infty$,
\begin{align*}
\Theta_N^-(t)&=\frac{4\,\pi\,\phi(N)^3}{\sqrt{|D_K|}\;t^3\,N^{3+\alpha}}
\sum_{p\in \OOO_K:\;0<|p|\leq N^\alpha,\;|p|<\frac{t}{2\lambda_\phi}N^\alpha} |p|^2
\\&=
  \frac{4\,\pi\,\phi(N)^3}{\sqrt{|D_K|}\;t^3\,N^{3+\alpha}}
  \Big(\frac{4\,\pi}{4\,\sqrt{|D_K|}}
  \min\Big\{N^\alpha,\frac{t\,N^\alpha}{2\lambda_\phi }\Big\}^4
  %\nonumber\\& \qquad\qquad\qquad\qquad\quad
  +\bigO\Big(\min\Big\{N^\alpha,\frac{t\,N^\alpha}{2\lambda_\phi }\Big\}^3
  \Big)\Big)\nonumber\\&= \frac{4\,\pi^2}{|D_K|}\Big(\frac{\phi(N)}{N^{1-\alpha}}
  \Big)^3\min\Big\{\frac{1}{t^3},\frac{t}{16\,\lambda_\phi^4}\Big\}
    +\bigO\Big(\frac{1}{N^\alpha}\Big(\frac{\phi(N)}{N^{1-\alpha}}
  \Big)^3\Big)\,.
\end{align*}
Since the function $t\mapsto \min\big\{\frac{1}{t^3}, \frac{t}{16\,
  \lambda_\phi^4}\big\}$ is bounded on $]0,+\infty[\,$, since
$\frac{\phi(N)}{N^{1-\alpha}}$ converges to $\lambda_\phi$ (hence is
bounded, uniformly in $t$) under the assumptions of Case
\hyperlink{cas2}{2}, and since we have $\big|\frac{1}{x}-\frac{1}{y}
\big| =\bigO(|x-y|)$ when $x,y$ remain in a compact subset of
$]0,+\infty[\,$, we have
\begin{align}
    \Theta_N^-(t)&=\frac{4\,\pi^2}{|D_K|}
  \min\Big\{\frac{\lambda_\phi^3}{t^3},\frac{t}{16\,\lambda_\phi}\Big\}
    +\bigO\big(\frac{1}{N^\alpha}\big)+\bigO\Big(
    \Big|\frac{\phi(N)}{N^{1-\alpha}}-\lambda_\phi\Big|\,\Big)
  \nonumber\\&=\frac{4\,\pi^2}{|D_K|}
  \min\Big\{\frac{\lambda_\phi^3}{t^3},\frac{t}{16\,\lambda_\phi}\Big\}
    +\bigO\big(\frac{1}{N^\alpha}\big)+\bigO\Big(
    \Big|\frac{N^{1-\alpha}}{\phi(N)}-\frac{1}{\lambda_\phi}\Big|\,\Big)\,.
    \label{eq:valueThetaNmoin}
\end{align}
Note that we have $\frac{\lambda_\phi^3}{t^3}\leq\frac{t}
{16\,\lambda_\phi}$ if and only if $t\geq 2\lambda_\phi$.

\medskip\noindent {\bf Step \hypertarget{step2}{2} of Case
  \hyperlink{cas2}{2}. } Let us now estimate $\Theta_N^+(t)$ as
$N\ra+\infty$ uniformly in $t\in[\varepsilon_N,+\infty[\,$.

\medskip
Let $t\in[\varepsilon_N,+\infty[$ and $N\in\NN\ssm\{0\}$. Note that
if $t> 2\lambda_\phi$, then $\Theta_N^+(t)=0$ (see Equation
\eqref{eq:defThetaNplu}). Hence we assume that $t\in [\varepsilon_N,
2\lambda_\phi]$ from now on in Step \hyperlink{step2}{2}. Let us
define two subsets in $\CC$ by
\begin{equation}\label{eq:defiANANprim}
\A_t=\Big\{z\in\CC:\frac{t}{2\lambda_\phi}\leq |z|\leq 1\Big\}
\quad
\text{and}
\quad
\A_{t,N}=\Big\{z\in\CC:\frac{t+\varepsilon_N}{2\lambda_\phi}
\leq|z|\leq 1\Big\}\,.
\end{equation}
The subset $\A_t$ (which is an annulus if $t<2\lambda_\phi$ and an
circle otherwise) contains $\A_{t,N}$ (which is an annulus if
$t<2\lambda_\phi- \varepsilon_N$, a circle if $t=2\lambda_\phi
-\varepsilon_N$, and is empty otherwise).

Note that the map $s\mapsto \frac{s}{s+\varepsilon_N}$ is increasing on
$[0,+\infty[\,$. Hence since ${\displaystyle \lim_{N\ra+\infty}}\,
\varepsilon_N=0$, for $N$ large enough, for every element $p\in\OOO_K$
occurring in the summation defining $\Theta_N^+(t)$ in Equation
\eqref{eq:defThetaNplu}, we have
\[
\frac{t\,N^\alpha}{2\,\lambda_\phi\,|p|}\leq \frac{t}{t+\varepsilon_N}\leq
\frac{2\lambda_\phi}{2\lambda_\phi+\varepsilon_N}
\leq 1- \frac{\varepsilon_N}{4\lambda_\phi}< 1\,.
\]
Hence for such elements $N$ and $p$, by the definition of
$\varepsilon_N$ in Equation \eqref{eq:defietaetap}, we have
\begin{align*}
\frac{N^{1-\alpha}}{\phi(N)}\,\frac{t\,N^\alpha}{2\,|p|}&\leq
\Big|\frac{N^{1-\alpha}}{\phi(N)}-\frac{1}{\lambda_\phi}\,\Big|
\,\frac{t\,N^\alpha}{2\,|p|}+
\frac{t\,N^\alpha}{2\,\lambda_\phi\,|p|}\leq
\Big|\frac{N^{1-\alpha}}{\phi(N)}-\frac{1}{\lambda_\phi}\,\Big|
\,\lambda_\phi+\frac{t\,N^\alpha}{2\,\lambda_\phi\,|p|}
\\&\leq \frac{\varepsilon_N}{8\lambda_\phi}+
1- \frac{\varepsilon_N}{4\lambda_\phi}=
1- \frac{\varepsilon_N}{8\lambda_\phi}< 1\,.
\end{align*}
Since $g'(u)=\bigO\big(\frac{1}{\sqrt{1-u}}\;\big)$ for $u\in [0,1[$
by Equation \eqref{eq:derivgvers1}, by the mean value theorem, and
since $\varepsilon_N\geq 8{\lambda_\phi}^2\, \big|\,
\frac{N^{1-\alpha}}{\phi(N)}-\frac{1}{\lambda_\phi}\,\big|$ by
Equation \eqref{eq:defietaetap}, for every $t\in
[\varepsilon_N, 2\lambda_\phi]$ and for every element $p\in\OOO_K$
occurring in the summation defining $\Theta_N^+(t)$, as $N\ra+\infty$,
we therefore have
\begin{align}
\Big|\;g\Big(\frac{N^{1-\alpha}}{\phi(N)}\,\frac{t\,N^\alpha}{2\,|p|}\Big)
-g\Big(\frac{t\,N^{\alpha}}{2\,\lambda_\phi\,|p|}\Big)\,\Big|
&=\bigO\Big(\frac{t\,N^{\alpha}}{|p|}\,
\Big|\,\frac{N^{1-\alpha}}{\phi(N)}-\frac{1}{\lambda_\phi}\,\Big|\,
\frac{1}{\sqrt{\varepsilon_N}}\,\Big)\nonumber\\&=
\bigO\Big(\frac{N^{\alpha}}{|p|}\,
\Big|\,\frac{N^{1-\alpha}}{\phi(N)}-\frac{1}{\lambda_\phi}\,\Big|^{\frac{1}{2}}
  \,\Big)\,.\label{eq:meanvaluegleq1}
\end{align}

Let us define
\begin{align}
\wt\Theta_N^+(t)&=\frac{1}{\sqrt{|D_K|}\;N^{\alpha}}
\sum_{p\in \OOO_K\cap (N^\alpha\A_{t,N})}\frac{1}{|p|}\,
g\Big(\frac{t\,N^\alpha}{2\,\lambda_\phi\,|p|}\Big)\,.
\label{eq:defiwtThetaNplu}
\end{align}
Since $t\geq \varepsilon_N$ in this Step \hyperlink{step2}{2} and
$\varepsilon_N\geq 8{\lambda_\phi}^2\, \big|\, \frac{N^{1-\alpha}}
{\phi(N)} \frac{1}{\lambda_\phi}\,\big|$, we have
\[
\frac{1}{t+\varepsilon_N}\leq \frac{1}{2\varepsilon_N}\leq
\frac{1}{16\,{\lambda_\phi}^2}\Big|\,
\frac{N^{1-\alpha}}{\phi(N)}-\frac{1}{\lambda_\phi}\,\Big|^{-1}\,.
\]
Hence by the definition of $\A_{t,N}$ in Equation
\eqref{eq:defiANANprim} and by Equation \eqref{eq:GaussSayousmoin2}
applied twice with $x=N^\alpha$ and with $x=\frac{t+\varepsilon_N}
{2\,\lambda_\phi}\,N^\alpha$, as $N\ra+\infty$, we have
\begin{align*}
  \sum_{p\in \OOO_K\cap (N^\alpha\A_{t,N})}\frac{1}{|p|^2}
  =\bigO\Big(\Big|\ln\frac{N^\alpha}{(t+\varepsilon_N)N^\alpha}\,\Big|\,\Big)
  +\bigO(1)=\bigO\Big(\Big|\ln\Big|\,\frac{N^{1-\alpha}}{\phi(N)}-
  \frac{1}{\lambda_\phi}\,\Big|\,\Big|\,\Big)\,.
\end{align*}
By Equations \eqref{eq:defThetaNplu}, \eqref{eq:defiwtThetaNplu} and
\eqref{eq:meanvaluegleq1}, for every $t\in[\varepsilon_N,
  2\lambda_\phi]$, as $N\ra+\infty$, we hence have
\begin{align}
  \big|\,\Theta_N^+(t)-\wt\Theta_N^+(t)\,\big|&=
  \bigO\Big(\frac{1}{N^{\alpha}}\sum_{p\in \OOO_K\cap (N^\alpha\A_{t,N})}
  \frac{1}{|p|}\,\frac{N^{\alpha}}{|p|}\,
  \Big|\,\frac{N^{1-\alpha}}{\phi(N)}-\frac{1}{\lambda_\phi}
  \,\Big|^{\frac{1}{2}}\,\Big)\nonumber\\&=\bigO\Big(\,
  \Big|\,\frac{N^{1-\alpha}}{\phi(N)}-\frac{1}{\lambda_\phi}\,
  \Big|^{\frac{1}{2}}\,\Big|\ln\Big|\,\frac{N^{1-\alpha}}{\phi(N)}-
  \frac{1}{\lambda_\phi}\,\Big|\,\Big|\,\Big)\,.
  \label{eq:ThetaNpluwtThetaNplu}
\end{align}

Let us consider the function $G_t:\A_t\ra[0,+\infty[$ defined by
\[
G_t:z\mapsto \frac{1}{|z|}\;
g\Big(\frac{t}{2\,\lambda_\phi\,|z|}\Big)\,.
\]
By the change of variables $z=N^{-\alpha}p$ in Equation
\eqref{eq:defiwtThetaNplu}, we have
\begin{align}
\wt\Theta_N^+(t)&=\frac{1}{\sqrt{|D_K|}\;N^{2\alpha}}
\sum_{z\in (N^{-\alpha}\OOO_K)\cap \A_{t,N}}G_t(z)\,.
\label{eq:valpartielThetaNplu}
\end{align}

Let us estimate the upper bound $\|d{G_t}\;_{\mid\A_{t,N}}\,\|_\infty$
on the operator norm of the (linear) differential $d_zG_t$ of $G_t$ at
every point $z\in \A_{t,N}$. Writing the elements of the annulus
$\A_t$ in polar coordinates $z=\rho \,e^{i\theta}$, the map $G_t$ does
not depend on the argument $\theta$. Since the function $\rho\mapsto
(1-\frac{t}{2\,\lambda_\phi \,\rho})^{-\frac{1}{2}}$ is nonincreasing
and since $\varepsilon_N\leq t$ in this Step \hyperlink{step2}{2}, for
every $\rho\in [\frac{t+\varepsilon_N} {2\lambda_\phi\,},1]$, we have
\[
\Big(1-\frac{t}{2\,\lambda_\phi\,\rho}\Big)^{-\frac{1}{2}} \leq
\Big(1-\frac{t}{2\,\lambda_\phi\,
  \frac{t+\varepsilon_N}{2\lambda_\phi\,}}\Big)^{-\frac{1}{2}} =
\Big(\frac{t+\varepsilon_N}{\varepsilon_N}\Big)^{\frac{1}{2}}
\leq\Big(\frac{2t}{\varepsilon_N}\Big)^{\frac{1}{2}}\,.
\]
Note that the map $s\mapsto\frac{s^{\frac{3}{2}}}{(s+
\varepsilon_N)^3}$ is nonincreasing on $[\varepsilon_N,+\infty[\,$,
hence is bounded by $\frac{1}{8}\varepsilon_N^{\;-\frac{3}{2}}$. Hence,
since the function $g$ is bounded and by Equation
\eqref{eq:derivgvers1}, as $N\ra+\infty$, for every $z=\rho\,
e^{i\theta}$ in the smaller subset $\A_{t,N}$ so that $\rho\geq
\frac{t+\varepsilon_N}{2\lambda_\phi\,}$, we have
\begin{align*}
\frac{\partial G_t}{\partial\rho}(z)&=
-\frac{1}{\rho^2}\,
g\Big(\frac{t}{2\,\lambda_\phi\,\rho}\Big)
-\frac{t}{2\,\lambda_\phi\,\rho^3}\,
g'\Big(\frac{t}{2\,\lambda_\phi\,\rho}\Big)\\&=
\bigO\Big(\frac{1}{(t+\varepsilon_N)^2}\Big)+
\bigO\Big(\frac{t}{(t+\varepsilon_N)^3}\,
\frac{1}{(1-\frac{t}{2\,\lambda_\phi\,\rho})^\frac{1}{2}}\Big)
\\&=\bigO\Big(\frac{1}{{\varepsilon_N}^2}\Big)+
\bigO\Big(\frac{t^{\frac{3}{2}}}{(t+\varepsilon_N)^3\;
  {\varepsilon_N}^\frac{1}{2}}\Big)=
\bigO\Big(\frac{1}{{\varepsilon_N}^2}\Big)\,.
\end{align*}
Therefore we have
\begin{equation}\label{eq:diffGN}
\big\|\,d{G_t}\,_{\mid \A_{t,N}}\,\big\|_\infty=
\Big\|\,\frac{\partial G_t}{\partial\rho}_{\mid \A_{t,N}}\,\Big\|_\infty
=\bigO\Big(\frac{1}{{\varepsilon_N}^2}\Big)\,.
\end{equation}
By Equation \eqref{eq:covoldiamsys}, the $\ZZ$-lattice 
$N^{-\alpha}\OOO_K$ of $\CC$ has covolume 
\begin{equation}\label{eq:covolNmoinsalphaOOOK}
\covol_{N^{-\alpha}\OOO_K}=
N^{-2\alpha} \covol_{\OOO_K}=\frac{\sqrt{|D_K|}}{2\,N^{2\alpha}}\,,
\end{equation}
and its diameter satisfies $\diam_{N^{-\alpha}\OOO_K}=
\bigO(N^{-\alpha})$. By the well-known approximation of integrals by
averages of sums over lattices points, we have
\begin{align*}
\int_{\A_{t,N}} G_t\;d\Leb_\CC&=
\covol_{N^{-\alpha}\OOO_K}\sum_{z\in (N^{-\alpha}\OOO_K)\cap \A_{t,N}} G_t(z)+
\bigO\big(\diam_{N^{-\alpha}\OOO_K} \|d{G_t}\,_{\mid \A_{t,N}}\|_\infty\big)\,.
\end{align*}
Hence by Equations \eqref{eq:covolNmoinsalphaOOOK} and
\eqref{eq:diffGN}, and since $\varepsilon_N\geq
\frac{1}{N^{\frac{\alpha}{4}}}$ by Equation \eqref{eq:defietaetap}, we
have
\begin{align}
  \sum_{z\in (N^{-\alpha}\OOO_K)\cap \A_{t,N}} G_t(z)&=
\frac{2\,N^{2\alpha}}{\sqrt{|D_K|}}\int_{\A_{t,N}} G_t\;d\Leb_\CC+
\bigO\Big(\frac{N^\alpha}{{\varepsilon_N}^2}\Big)
\nonumber\\&\label{eq:sumintGN}
=\frac{2\,N^{2\alpha}}{\sqrt{|D_K|}}\int_{\A_{t,N}} G_t\;d\Leb_\CC+
\bigO(N^{\frac{3\alpha}{2}})\,.
\end{align}

By the definition of the function $g$ in Equation \eqref{eq:defig} and
by an elementary computation, a primitive of the map $s\mapsto
\frac{g(s)}{s^2}$ on $]0,1]$ is the map $h:\;]0,1]\ra\RR$ defined by
\begin{equation}\label{eq:defih}
\forall\;u\in\;]0,1],\quad h(u)=
    \frac{1}{4u^4}\big(u(1-2u^2)\sqrt{1-u^2}-\arcsin(u)\big)\,.
\end{equation}
Furthermore, we have
\[
h(1)=-\frac{\pi}{8}\qquad\text{and}\qquad
h(u)=-\frac{2}{3\,u} +\bigO(u)\text{ as }u\ra 0\,.
\]
Since $d\Leb_\CC(z)=\rho\;d\rho\,d\theta$ and by the change of
variable $s=\frac{t}{2\,\lambda_\phi\,\rho}$, we have
\begin{align}
  \int_{\A_t} G_t\;d\Leb_\CC
  &=2\,\pi\int_{\frac{t}{2\lambda_\phi}\leq\rho \leq 1}
  \frac{1}{\rho} \,g\Big(\frac{t}{2\,\lambda_\phi\,\rho}
  \Big)\,\rho\;d\rho=\frac{\pi\,t}{\lambda_\phi}
  \int_{\frac{t}{2\lambda_\phi}\leq s \leq1}\frac{g(s)}{s^2}\;ds
  \nonumber\\&=\frac{\pi\,t}{\lambda_\phi}\Big(h(1)-
  h\Big(\frac{t}{2\lambda_\phi}\Big)\Big)\,.\label{eq:caclintGt}
\end{align}

Since the function $s\mapsto sg(s)$ is increasing on $[0,1]$, the
function $G_t$ on $\A_t$ is uniformly bounded by $\frac{2\lambda_\phi}
{t} g(1)=\frac{\pi\,\lambda_\phi}{t}$. Hence, by Equation
\eqref{eq:defiANANprim} and since $t\geq \varepsilon_N$ in this Step
\hyperlink{step2}{2}, when $N\ra +\infty$, we have
\begin{align}
  \int_{\A_t\smallsetminus \A_{t,N}} G_t\;d\Leb_\CC&=
  \bigO\big(\frac{1}{t}\operatorname{Area}(\A_t\ssm\A_{t,N})\big)
  =\bigO\Big(\frac{\pi}{t}\Big(\frac{(t+\varepsilon_N)^2}{4{\lambda_\phi}^{2}}
  -\frac{t^2}{4{\lambda_\phi}^{2}}\Big)\Big)\nonumber\\&=
  \bigO\big(\varepsilon_N+\frac{\varepsilon_N^2}{t}\big)=\bigO(\varepsilon_N)
  \,.\label{eq:intpetitann}
\end{align}
Therefore by Equations \eqref{eq:sumintGN}, \eqref{eq:caclintGt} and
\eqref{eq:intpetitann}, we have
\begin{align*}
  \sum_{z\in (N^{-\alpha}\OOO_K)\cap\A_{t,N}}G_t(z)&=
  \sum_{z\in (N^{-\alpha}\OOO_K)\cap\A_{t,N}}G_t(z)-\frac{2\,N^{2\alpha}}{\sqrt{|D_K|}}
  \int_{\A_{t,N}} G_t\;d\Leb_\CC\\&\quad+\frac{2\,N^{2\alpha}}{\sqrt{|D_K|}}
  \int_{\A_t}G_t\;d\Leb_\CC -\frac{2\,N^{2\alpha}}{\sqrt{|D_K|}}
  \int_{\A_t\smallsetminus \A_{t,N}}G_t\;d\Leb_\CC
  \\&=\bigO(N^{\frac{3\alpha}{2}})+
  \frac{2\,N^{2\alpha}\,\pi\,t}{\sqrt{|D_K|}\;\lambda_\phi}\Big(h(1)-
  h\Big(\frac{t}{2\lambda_\phi}\Big)\Big)
  +\bigO(N^{2\alpha}\,\varepsilon_N)\,.
\end{align*}
Thus by Equations \eqref{eq:ThetaNpluwtThetaNplu} and
\eqref{eq:valpartielThetaNplu}, by the definition of $\varepsilon_N$
in Equation \eqref{eq:defietaetap}, as $N\ra+\infty$, for every
$t\in[\varepsilon_N,2\lambda_\phi]$, we have
\begin{align}
  &\Theta_N^+(t)=\wt\Theta_N^+(t)+\bigO\Big(\,
  \Big|\frac{N^{1-\alpha}}{\phi(N)}-\frac{1}{\lambda_\phi}\Big|^{\frac{1}{2}}
  \,\Big|\ln\Big|\,\frac{N^{1-\alpha}}{\phi(N)}-\frac{1}{\lambda_\phi}\Big|
  \,\Big|\,\Big)\nonumber\\=\;&\frac{1}{\sqrt{|D_K|}\;N^{2\alpha}}
\sum_{z\in (N^{-\alpha}\OOO_K)\cap \A_{t,N}}G_t(z)+\bigO\Big(\,
\Big|\frac{N^{1-\alpha}}{\phi(N)}-\frac{1}{\lambda_\phi}\Big|^{\frac{1}{2}}
\,\Big|\ln\Big|\,\frac{N^{1-\alpha}}{\phi(N)}-\frac{1}{\lambda_\phi}\Big|
  \,\Big|\,\Big)\nonumber\\=\;&\frac{2\,\pi\,t}{|D_K|\,\lambda_\phi}
  \Big(h(1)-h\Big(\frac{t}{2\lambda_\phi}\Big)\Big)+\bigO(\varepsilon_N)
  +\bigO\big(\frac{1}{N^{\frac{\alpha}{2}}}\big)+\bigO\Big(\,
  \Big|\frac{N^{1-\alpha}}{\phi(N)}-\frac{1}{\lambda_\phi}\Big|^{\frac{1}{2}}
  \Big|\ln\Big|\,\frac{N^{1-\alpha}}{\phi(N)}-\frac{1}{\lambda_\phi}\,\Big|\,
  \Big|\,\Big)\nonumber\\=\;&\frac{2\,\pi\,t}{|D_K|\,\lambda_\phi}
  \Big(h(1)-h\Big(\frac{t}{2\lambda_\phi}\Big)\Big)+
  \bigO\big(\frac{1}{N^{\frac{\alpha}{4}}}\big)+\bigO\Big(\,
  \Big|\frac{N^{1-\alpha}}{\phi(N)}-\frac{1}{\lambda_\phi}\Big|^{\frac{1}{2}}
  \Big|\ln\Big|\,\frac{N^{1-\alpha}}{\phi(N)}-
  \frac{1}{\lambda_\phi}\,\Big|\,\Big|\,\Big)\,.
  \label{eq:valueThetaNplu}
\end{align}

\medskip\noindent {\bf Step \hypertarget{step3}{3} of Case
  \hyperlink{cas2}{2}. } Let us finally estimate $\Theta_N^0(t)$ as
$N\ra+\infty$ uniformly in $t\in[0,+\infty[\,$.

\medskip
Let $N\in\NN\ssm\{0\}$. Note that if $t>2\lambda_\phi$, then
$\Theta_N^0(t)=0$ (see Equation \eqref{eq:defThetaNzero}). Hence we
assume from now on in Step \hyperlink{step3}{3} that $t\in[0,
  2\lambda_\phi]$. Using the notation of Equation
\eqref{eq:defiANANprim}, Equation \eqref{eq:defThetaNzero} may be
written
\begin{align*}
\Theta_N^0(t)&=\frac{1}{\sqrt{|D_K|}\;N^{\alpha}}
\sum_{p\in \OOO_K\cap (N^{\alpha}\A_t\smallsetminus N^{\alpha}\A_{t,N})}
\frac{1}{|p|}\,\,g\Big(\frac{t\,N}{2\,\phi(N)\,|p|}\Big)\,.
\end{align*}
Since the nonnegative function $g$ is uniformly bounded from above (by
$\frac{\pi}{2}$), using twice Equation \eqref{eq:GaussSayous} with
$\beta'=-1$ and $x=\frac{(t+\varepsilon_N)}{2\lambda_\phi}
N^{\alpha}$, $x=\frac{t}{2\lambda_\phi}N^{\alpha}$, as $N\ra+\infty$,
we have
\begin{align}
  \Theta_N^0(t)&=\bigO\Big(\frac{1}{N^{\alpha}}\sum_{p\in \OOO_K\cap
    (N^{\alpha}\A_t\smallsetminus N^{\alpha}\A_{t,N})} \frac{1}{|p|}\Big)
  \nonumber\\&=\bigO\Big(\frac{1}{N^{\alpha}}
  \Big(\frac{(t+\varepsilon_N)}{2\lambda_\phi}N^{\alpha}-
  \frac{t}{2\lambda_\phi}N^{\alpha}\Big)\Big)=
  \bigO\big(\varepsilon_N\big)\,.\label{eq:valueThetaNzero}
\end{align}

\medskip\noindent{\bf Conclusion of Case \hyperlink{cas2}{2}. } To
conclude, we gather the estimates of Steps \hyperlink{step1}{1},
\hyperlink{step2}{2} and \hyperlink{step3}{3}.

Let us first compute $\Theta_N(t)$ for $t\in\; [\varepsilon_N, +
\infty[\,$. We separate the computation into the case $t>
2\lambda_\phi$ and the case $\varepsilon_N\leq t\leq 2\lambda_\phi$.

If $t> 2\lambda_\phi$, then $\Theta_N^0(t) =\Theta_N^+(t)=0$ as seen in
Steps \hyperlink{step2}{2} and \hyperlink{step3}{3}. By Equations
\eqref{eq:ThetaNsumpluzeromoi} and \eqref{eq:valueThetaNmoin}, and by
the value of $\rho_{1-\alpha}(t)$ when $t> 2\lambda_\phi$ given in
Equation \eqref{eq:defpaircorrfunctgene}, we hence have
\begin{align}
  \Theta_N(t)&=\Theta_N^-(t)=\frac{4\,\pi^2\,\lambda_\phi^3}{|D_K|\,t^3}
    +\bigO\big(\frac{1}{N^\alpha}\big)+\bigO\Big(\,
    \Big|\frac{N^{1-\alpha}}{\phi(N)}-\frac{1}{\lambda_\phi}\Big|\,\Big)
    \nonumber\\&=\rho_{1-\alpha}(t)+
    \bigO\big(\frac{1}{N^{\frac{\alpha}{4}}}\big)+\bigO\Big(\,
  \Big|\frac{N^{1-\alpha}}{\phi(N)}-\frac{1}{\lambda_\phi}\Big|^{\frac{1}{2}}
  \Big|\ln\Big|\,\frac{N^{1-\alpha}}{\phi(N)}-\frac{1}{\lambda_\phi}
  \,\Big|\,\Big|\,\Big)\,.\label{eq:valueThetaNfinalbigt}
\end{align}
      
If $\varepsilon_N\leq t\leq 2\lambda_\phi$, by plugging in Equation
\eqref{eq:ThetaNsumpluzeromoi} the three Equations
\eqref{eq:valueThetaNmoin}, \eqref{eq:valueThetaNzero} and
\eqref{eq:valueThetaNplu}, by the definition of $\varepsilon_N$ in
Equation \eqref{eq:defietaetap}, by the value of the function $h$
given in Equation \eqref{eq:defih} (with $h(1)=-\frac{\pi}{8}$) and by
the value of $\rho_{1-\alpha}(t)$ when $t\leq 2\lambda_\phi$ given in
Equation \eqref{eq:defpaircorrfunctgene}, as $N\ra+\infty$, we have
\begin{align}
  &\Theta_N(t)=\Theta_N^-(t)+\Theta_N^0(t)+\Theta_N^+(t)\nonumber\\=\;&
  \left(\frac{\pi^2\,t}{4\,|D_K|\,\lambda_\phi}
    +\bigO\big(\frac{1}{N^\alpha}\big)+\bigO\Big(\,
    \Big|\frac{N^{1-\alpha}}{\phi(N)}-\frac{1}{\lambda_\phi}\Big|\,\Big)\right)
    +\bigO\big(\varepsilon_N\big)\nonumber\\&+\left(
    \frac{2\,\pi\,t}{|D_K|\,\lambda_\phi}
  \Big(h(1)-h\Big(\frac{t}{2\lambda_\phi}\Big)\Big)+
  \bigO\big(\frac{1}{N^{\frac{\alpha}{4}}}\big)+\bigO\Big(\,
  \Big|\frac{N^{1-\alpha}}{\phi(N)}-\frac{1}{\lambda_\phi}\Big|^{\frac{1}{2}}
  \Big|\ln\Big|\,\frac{N^{1-\alpha}}{\phi(N)}-
  \frac{1}{\lambda_\phi}\,\Big|\,\Big|\,\Big)\right)\nonumber\\=\;&
   \frac{8\,\pi\,\lambda_\phi^3}{|D_K|\,t^3}
\Big(\arcsin\Big(\frac{t}{2\lambda_\phi}\Big)-
\frac{t}{2\lambda_\phi}\Big(1-\frac{t^2}{2\lambda_\phi^2}\Big)
\Big(1-\frac{t^2}{4\lambda_\phi^2}\Big)^{\frac{1}{2}}\,\Big)\nonumber\\&+
  \bigO\big(\frac{1}{N^{\frac{\alpha}{4}}}\big)+\bigO\Big(\,
  \Big|\frac{N^{1-\alpha}}{\phi(N)}-\frac{1}{\lambda_\phi}\Big|^{\frac{1}{2}}
  \Big|\ln\Big|\,\frac{N^{1-\alpha}}{\phi(N)}-\frac{1}{\lambda_\phi}\,\Big|
  \,\Big|\,\Big)\nonumber\\=\;& \rho_{1-\alpha}(t)+
   \bigO\big(\frac{1}{N^{\frac{\alpha}{4}}}\big)+\bigO\Big(\,
  \Big|\frac{N^{1-\alpha}}{\phi(N)}-\frac{1}{\lambda_\phi}\Big|^{\frac{1}{2}}
  \Big|\ln\Big|\,\frac{N^{1-\alpha}}{\phi(N)}-\frac{1}{\lambda_\phi}\,\Big|
  \,\Big|\,\Big)\,.\label{eq:valueThetaNfinal}
\end{align}
By Equation \eqref{eq:nettoyageThetaNcas2}, since the function $g$ is
bounded on $[0,+\infty[$ and by Equation \eqref{eq:GaussSayous} with
$\beta'=-1$ and $x=N^\alpha$, as $N\ra+\infty$, uniformly in $t\in
[0,+\infty[\,$, we have
\[
\Theta_N(t)=\bigO\Big(\frac{1}{N^{\alpha}}
\sum_{p\in \OOO_K:\;0<|p|\leq N^\alpha} \frac{1}{|p|}\Big)= \bigO(1)\,.
\]
We have already seen after Equation \eqref{eq:defpaircorrfunctgene}
that the function $\rho_{1-\alpha}$ is bounded on $[0,+\infty[\,$,
hence $\int_0^{\varepsilon_N} \rho_{1-\alpha}(t)\,f(t)\;dt=\bigO(
\varepsilon_N \,\|f\|_\infty)$.  Therefore, since the support of $f$ is
contained in $[0,A]$, by Equation \eqref{eq:valueThetaNfinal}, and by
the definition of $\varepsilon_N$ in Equation \eqref{eq:defietaetap},
we have
\begin{align}
  \int_0^{+\infty} \Theta_N(t)\,f(t)\;dt&=\int_0^{\varepsilon_N}
  \Theta_N(t)\,f(t)\;dt+\int_{\varepsilon_N}^{A}
  \Theta_N(t)\,f(t)\;dt \nonumber\\&=
  \bigO(\varepsilon_N\,\|f\|_\infty)+\int_{\varepsilon_N}^{A}
  \rho_{1-\alpha}(t)\,f(t)\;dt\nonumber\\&\quad+
  \bigO\Big(\Big(\frac{1}{N^{\frac{\alpha}{4}}}+
  \Big|\frac{N^{1-\alpha}}{\phi(N)}-\frac{1}{\lambda_\phi}\Big|^{\frac{1}{2}}
  \Big|\ln\Big|\,\frac{N^{1-\alpha}}{\phi(N)}-\frac{1}{\lambda_\phi}\,
  \Big|\,\Big|\,\Big)\int_{\varepsilon_N}^{A}
  |f(t)|\;dt\Big)\nonumber\\&= \int_0^{+\infty}
  \rho_{1-\alpha}(t)\,f(t)\;dt\nonumber\\&\quad+
  \bigO\Big(\Big(\frac{1}{N^{\frac{\alpha}{4}}}+
  \Big|\frac{N^{1-\alpha}}{\phi(N)}-\frac{1}{\lambda_\phi}\Big|^{\frac{1}{2}}
  \Big|\ln\Big|\,\frac{N^{1-\alpha}}{\phi(N)}-\frac{1}{\lambda_\phi}\,
  \Big|\,\Big|\,\Big)A\,\|f\|_\infty\Big)\,.\label{eq:ThetaNcas2integ}
\end{align}

As $N\ra+\infty$, since $\psi(N)=\frac{N^{3+\alpha}}{\phi(N)}$, since
we have $\frac{\phi(N)}{N^{1-\alpha}} =\bigO(1)$ under the assumptions
of Case \hyperlink{cas2}{2}, and since $\ga< \frac{1-2\alpha}{2}$, we
have
\begin{align}
  \frac{\phi(N)\,N^{4\alpha+2\ga}}{\psi(N)}=
  \frac{\phi(N)^2}{N^{3-3\alpha-2\ga}}=
  \bigO\Big(\frac{1}{N^{1-\alpha-2\ga}}\Big)=
  \bigO\Big(\frac{1}{N^{\alpha}}\Big)\,.
  \label{eq:controlerror1cas2}
\end{align}
Similarly, we have
\begin{equation}\label{eq:controlerror2cas2}
  \frac{N^{2+2\alpha-2\ga}}{\psi(N)}=\frac{\phi(N)}{N^{1-\alpha+2\ga}}=
  \bigO\Big(\frac{1}{N^{2\ga}}\Big)\,.
\end{equation}
Therefore, using in Equation \eqref{eq:caspasplusinfty} the three Equations
\eqref{eq:ThetaNcas2integ}, \eqref{eq:controlerror1cas2} and
\eqref{eq:controlerror2cas2}, we have
\begin{align}
  \R^+_{N}(f)&=\int_0^{+\infty} \rho_{1-\alpha}(t)\,f(t)\;dt
  \quad+
  \bigO\Big(\frac{1}{N^{\alpha}}\,\|f'\|_\infty\Big)\nonumber\\&+
  \bigO\Big(\max\Big\{\frac{A}{N^{\frac{\alpha}{4}}},\;A\,
  \Big|\frac{N^{1-\alpha}}{\phi(N)}-\frac{1}{\lambda_\phi}\Big|^{\frac{1}{2}}
  \Big|\ln\Big|\,\frac{N^{1-\alpha}}{\phi(N)}-\frac{1}{\lambda_\phi}\,
  \Big|\,\Big|,\;\frac{1}{N^{2\ga}}\Big\}\,\|f\|_\infty\Big)\,.
\label{eq:cas2plus}
\end{align}
By symmetry, for every $f\in C_c(\RR)$, this proves the second case of
Theorem \ref{theo:sumsquaremain}.

When $\phi:N\ra N^\beta$ is a power function, the assumptions of Case
\hyperlink{cas2}{2} on $\phi$ are satisfied if and only if
$\beta=1-\alpha$ and $D_K\equiv 0\!\!\mod 4$, and we have
$\lambda_\phi=1$, so that the value of $\rho_{1-\alpha}$ given by
Equation \eqref{eq:defpaircorrfunctgene} becomes the one given in
Equation \eqref{eq:defpaircorrfunct} of the Introduction. This proves
the second case of Theorem \ref{theo:mainintro}.

\medskip
\noindent{\bf Case \hypertarget{cas3}{3}. } Let us assume that
${\displaystyle\lim_{N\ra+\infty}} \;\frac{N^{1-\alpha}}{\phi(N)} =0$,
that ${\displaystyle\limsup_{N\ra+\infty}}\; \frac{\phi(N)}
{N^{1-\frac{\alpha}{2}}} <+\infty$ and that $\alpha<\frac{2}{5}$.

Let us take $\psi$ to be the function $N\mapsto \frac{N^3\;S_{N,0}}
{\phi(N)}$, with $S_{N,0}$ defined in Equation \eqref{eq:defiSN} for
$k=0$. Let $\ga =\frac{2-3\alpha}{8}$. The assumption
\eqref{eq:conditiongaeps} is then satisfied since $\alpha <
\frac{2}{5}$.  Since $\frac{\phi(N)}{N^{1-\frac{\alpha}{2}}}
=\bigO(1)$ under the assumptions of Case \hyperlink{cas3}{3}, and
since $\alpha < \frac{2}{3}$, we have
\[
\frac{\phi(N)}{N^{2-2\alpha-2\ga}}=
\bigO\Big(\frac{1}{N^{1-\frac{3}{2}\alpha-2\ga}}\Big)=
\bigO\Big(\frac{1}{N^{\frac{2-3\alpha}{4}}}\Big)\,,
\]
which tends to $0$ as $N\ra+\infty$. Hence the assumption
\eqref{eq:scalingrange} is satisfied and we may apply the results of
Section \ref{sec:linear}.

We are going to prove that in Case \hyperlink{cas3}{3}, as $N\ra
+\infty$, the function $\Theta_N$ defined in Equation
\eqref{eq:defiThetaN} converges uniformly on compact subsets of
$[0,+\infty[$ to the constant function $\frac{8\,\pi}{3\,|D_K|}$.
Below are the graphs of $\Theta_N$ for various $N$, for the power
scaling $\phi(N)=N^\beta$, in the Gaussian case $K=\QQ[i]$, with
$\alpha=0.15$ and $\beta =0.9\in\;]1-\alpha,1-\frac{\alpha}{2}]$, and
in dashed black the horizontal line with height $\frac{8\,\pi}{3|D_K|}
\simeq 2.094$.
\begin{center}
    \includegraphics[width=15cm]{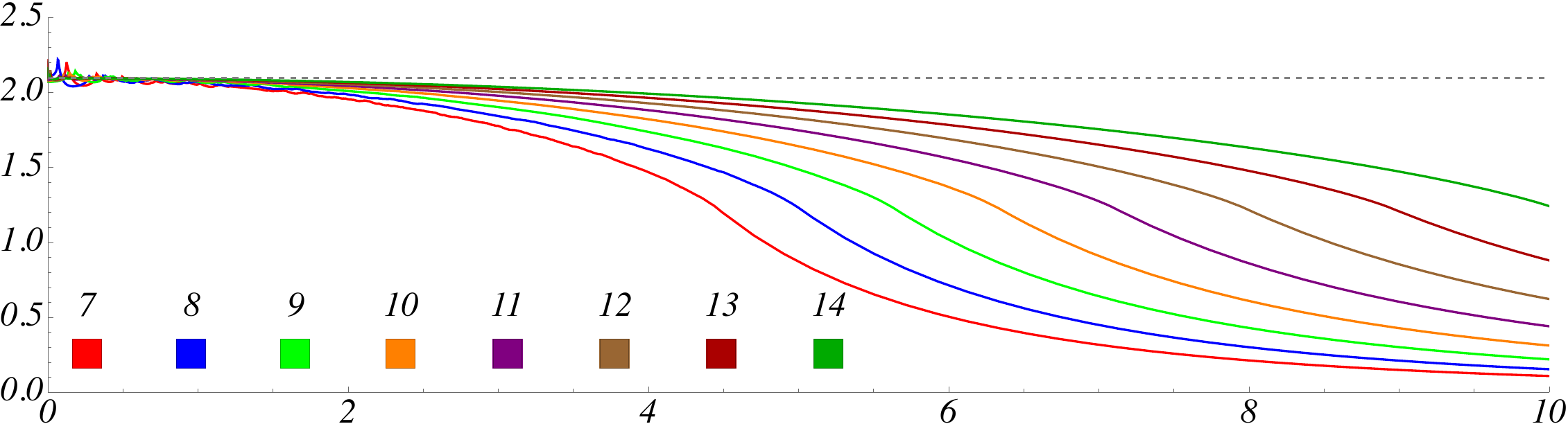}
  \end{center}\vspace*{-0.4cm} \begin{center}
\addtocounter{fig}{1} {\small Figure \arabic{fig}~: Graph of
  $\Theta_N$ for $N=10^m$ where $m$ is $7$ (red), $8$ (blue) up to
  $14$ (dark green).}
\end{center}

For every $N\in\NN\ssm\{0\}$, let us define
\begin{equation}\label{eq:defiMsubN}
M_N=N^\alpha\Big(\,\frac{N^{1-\alpha}}{\phi(N)}\,\Big)^{\frac{2}{3}}\,.
\end{equation}
Since $M_N=\big(\frac{N^{1-\frac{\alpha}{2}}}{\phi(N)}
\big)^{\frac{2}{3}} N^{\frac{2\alpha}{3}}$ and since
$\frac{N^{1-\frac{\alpha}{2}}}{\phi(N)}$ is bounded from below by a
positive constant under the assumptions of Case \hyperlink{cas3}{3},
we have ${\displaystyle\lim_{N\ra+\infty}}\,M_N =+\infty$. Let
\[
S(M_N)=\frac{|D_K|}{4\,\pi}\sum_{p\in \OOO_K:\;0<|p|\leq M_N}
\frac{1}{c'_p\,|v_p|}\,.
\]
As in order to obtain Equation \eqref{eq:controlSN} with $k=0$, we
have
\[
S(M_N)=\bigO(M_N)\,.
\]
Let $t\in[0,A]$. Since $\psi(N)=\frac{N^3\;S_{N,0}} {\phi(N)}$, since
the function $g$ is bounded, by Equation \eqref{eq:controlSN} with $k=
0$, and by the definition of $M_N$ in Equation \eqref{eq:defiMsubN},
we hence have
\begin{align}
\frac{N^3}{\psi(N)\,\phi(N)}
\sum_{p\in \OOO_K:\;0<|p|\leq M_N} \frac{1}{c'_p\,|v_p|}
\,g\Big(\frac{t\,N}{2\,\phi(N)\,|p|}\Big)
&=\bigO\Big(\frac{S(M_N)}{S_{N,0}}\Big)
=\bigO\Big(\frac{M_N}{N^{\alpha}}\Big)\nonumber
\\&=\bigO\Big(\Big(\frac{N^{1-\alpha}}{\phi(N)}\Big)^\frac{2}{3}\,\Big)\,.
\label{eq:ThetaNmoincas3}
\end{align}
Note that if $p\in\OOO_K$ satisfies $|p|\geq M_N$, then we have
$\frac{t\,N}{2\,\phi(N)\,|p|}\leq\frac{A\,N}{2\,\phi(N)\,M_N}$. Since
$\psi(N)=\frac{N^3\;S_{N,0}}{\phi(N)}$, since we have $g(u)=
\frac{2}{3} +\bigO(u^2)$ for $u\in[0,+\infty[\,$, by the definition of
$M_N$ in Equation \eqref{eq:defiMsubN} and since $A\geq 1$, we
hence have
\begin{align}
  &\frac{N^3}{\psi(N)\,\phi(N)}
  \sum_{p\in \OOO_K:\;M_N<|p|\leq N^{\alpha}} \frac{1}{c'_p\,|v_p|}
\,g\Big(\frac{t\,N}{2\,\phi(N)\,|p|}\Big)\nonumber
\\=\;&\frac{8\,\pi}{3\,|D_K|}\frac{S_{N,0}-S(M_N)}{S_{N,0}}+
\bigO\Big(\Big(\frac{A\,N}{\phi(N)\,M_N}\Big)^2\,\Big)\nonumber\\=\;&
\frac{8\,\pi}{3\,|D_K|}+\bigO\Big(\frac{M_N}{N^{\alpha}}\Big)+
\bigO\Big(A^2\Big(\frac{N^{1-\alpha}}{\phi(N)}
\frac{N^\alpha}{M_N}\Big)^2\,\Big)\nonumber
\\=\;&\frac{8\,\pi}{3\,|D_K|}+
\bigO\Big(A^2\Big(\frac{N^{1-\alpha}}{\phi(N)}\Big)^\frac{2}{3}\,\Big)\,.
\label{eq:ThetaNplucas3}
\end{align}
By decomposing the index set of the sum defining the function
$\Theta_N$ as
\[
\{p\in\OOO_K:\;0<|p|\leq N^\alpha\} =\{p\in\OOO_K:\;0<|p|\leq M_N\}
\sqcup\{p\in\OOO_K:\;M_N<|p|\leq N^{\alpha}\}\,,
\]
it follows from Equations \eqref{eq:defiThetaN},
\eqref{eq:ThetaNmoincas3} and \eqref{eq:ThetaNplucas3} that as
$N\ra+\infty$, for every $t\in [0,A]$, we have
\[
\Theta_N(t)=\frac{8\,\pi}{3\,|D_K|}+
\bigO\Big(A^2\Big(\frac{N^{1-\alpha}}{\phi(N)}\Big)^\frac{2}{3}\,\Big)\,.
\]
Therefore, since the support of $f$ is contained in $[0,A]$, we have
\begin{equation}\label{eq:ThetaNcas3integ}
  \int_0^{+\infty} \Theta_N(t)\,f(t)\;dt=
  \int_0^{+\infty} \frac{8\,\pi}{3\,|D_K|}\,f(t)\;dt+
  \bigO\Big(A^3\Big(\frac{N^{1-\alpha}}{\phi(N)}\Big)^\frac{2}{3}
  \|f\|_\infty\,\Big)\,.
\end{equation}

As $N\ra+\infty$, since $\psi(N)=\frac{N^3\;S_{N,0}} {\phi(N)}$, by
Equation \eqref{eq:controlSN} with $k=0$, since we have
$\frac{\phi(N)}{N^{1-\frac{\alpha}{2}}} =\bigO(1)$ under the
assumptions of Case \hyperlink{cas3}{3}, and since $\ga=
\frac{2-3\alpha}{8}$, we have
\begin{align}
  \frac{\phi(N)\,N^{4\alpha+2\ga}}{\psi(N)}&=
  \frac{\phi(N)^2\,N^{4\alpha+2\ga}}{N^3\,S_{N,0}}=
  \bigO\Big(\frac{\phi(N)^2}{N^{3-3\alpha-2\ga}}\Big)\nonumber\\&=
  \bigO\Big(\frac{1}{N^{1-2\alpha-2\ga}}\Big)=
  \bigO\Big(\frac{1}{N^{\frac{2-5\alpha}{4}}}\Big)\,.
  \label{eq:controlerror1cas3}
\end{align}
Similarly, we have
\begin{equation}\label{eq:controlerror2cas3}
  \frac{N^{2+2\alpha-2\ga}}{\psi(N)}=
  \frac{N^{2+2\alpha-2\ga}\,\phi(N)}{N^3\,S_{N,0}}=
  \bigO\Big(\frac{\phi(N)}{N^{1-\alpha+2\ga}}\Big)=
  \bigO\Big(\frac{1}{N^{2\ga-\frac{\alpha}{2}}}\Big)=
  \bigO\Big(\frac{1}{N^{\frac{2-5\alpha}{4}}}\Big)\,.
\end{equation}
Therefore, using in Equation \eqref{eq:caspasplusinfty} the three Equations
\eqref{eq:ThetaNcas3integ}, \eqref{eq:controlerror1cas3} and
\eqref{eq:controlerror2cas3}, we have
\begin{align}
  \R^+_{N}(f)&=\int_0^{+\infty} \frac{8\,\pi}{3\,|D_K|}\,f(t)\;dt
  \nonumber\\&\quad+
  \bigO\Big(\frac{1}{N^{\frac{2-5\alpha}{4}}}\,\|f'\|_\infty\Big)+
  \bigO\Big(\max\Big\{A^3\Big(\frac{N^{1-\alpha}}{\phi(N)}\Big)^\frac{2}{3}
  ,\;\frac{1}{N^{\frac{2-5\alpha}{4}}}\Big\}\,\|f\|_\infty\Big)\,.
\label{eq:cas3plus}
\end{align}
By symmetry, for every $f\in C_c(\RR)$, this proves the third case of
Theorem \ref{theo:sumsquaremain}. 

When $\phi:N\ra N^\beta$ is a power function, the assumptions of Case
\hyperlink{cas3}{3} on $\phi$ are satisfied if and only if $\beta\in\;
  ]1-\alpha, 1-\frac{\alpha}{2}]$. If furthermore $D_K\equiv 0\!\!\mod
4$, then by Equation \eqref{eq:defiSNcasdiscnul} with $k=0$, as
$N\ra+\infty$, we have $\psi(N)=\frac{N^3\,S_{N,0}}{\phi(N)}\sim
N^{3+\alpha-\beta}$. This proves the third case of Theorem
\ref{theo:mainintro} when $\beta\in\;]1-\alpha, 1-\frac{\alpha}{2}]$.

\medskip
\noindent{\bf Case \hypertarget{cas5}{5}. } Let us assume that
${\displaystyle\limsup_{N\ra+\infty}} \frac{N}{\phi(N)} <+\infty\,$,
that $\lambda'''_N={\displaystyle \lim_{N\ra+\infty}}\;\frac{\phi(N)}
{N^{1+\frac{\alpha}{2}}}=0$ and that $\alpha\leq\frac{2}{11}$.

Let us take $\psi:N\mapsto \frac{N^3\,S_{N,0}}{\phi(N)}$ and
$\ga=\frac{1-4\alpha}{2}$. The assumptions \eqref{eq:conditiongaeps}
and \eqref{eq:scalingrange} are satisfied (since $\frac{\phi(N)}
{N^{2-2\alpha-2\ga}}= \frac{\phi(N)} {N^{1+\frac{\alpha}{2}}}\,
N^{-\frac{3\alpha}{2}}$ tends to $0$ as $N\ra+\infty$ under the
assumptions of Case \hyperlink{cas5}{5}), hence we may apply the
results of Section \ref{sec:linear}. For every $N\in\NN\ssm\{0\}$, let
\[
\eta_N=\frac{N}{\phi(N)}\,,
\]
and note that $\eta_N$ remains bounded as $N\ra+\infty$ under the
assumptions of Case \hyperlink{cas5}{5}.

We are going to prove that in this case, the function $\Theta_N$
defined in Equation \eqref{eq:defiThetaN} converges uniformly on
compact subsets of $[0,+\infty[$ to the constant function
$\frac{8\pi}{3|D_K|}$. Below are the graphs of $\Theta_N$ for various
$N$, for the power scaling $\phi(N)=N^\beta$, in the Gaussian case
$K=\QQ[i]$, with $\alpha=0.15$ and $\beta =1.05\in[1,1+
\frac{\alpha}{2}[\,$.
%, and in black the horizontal line with height
%$\frac{8\pi}{3|D_K|} \simeq 2.094$.
\begin{center}
    \includegraphics[width=15cm]{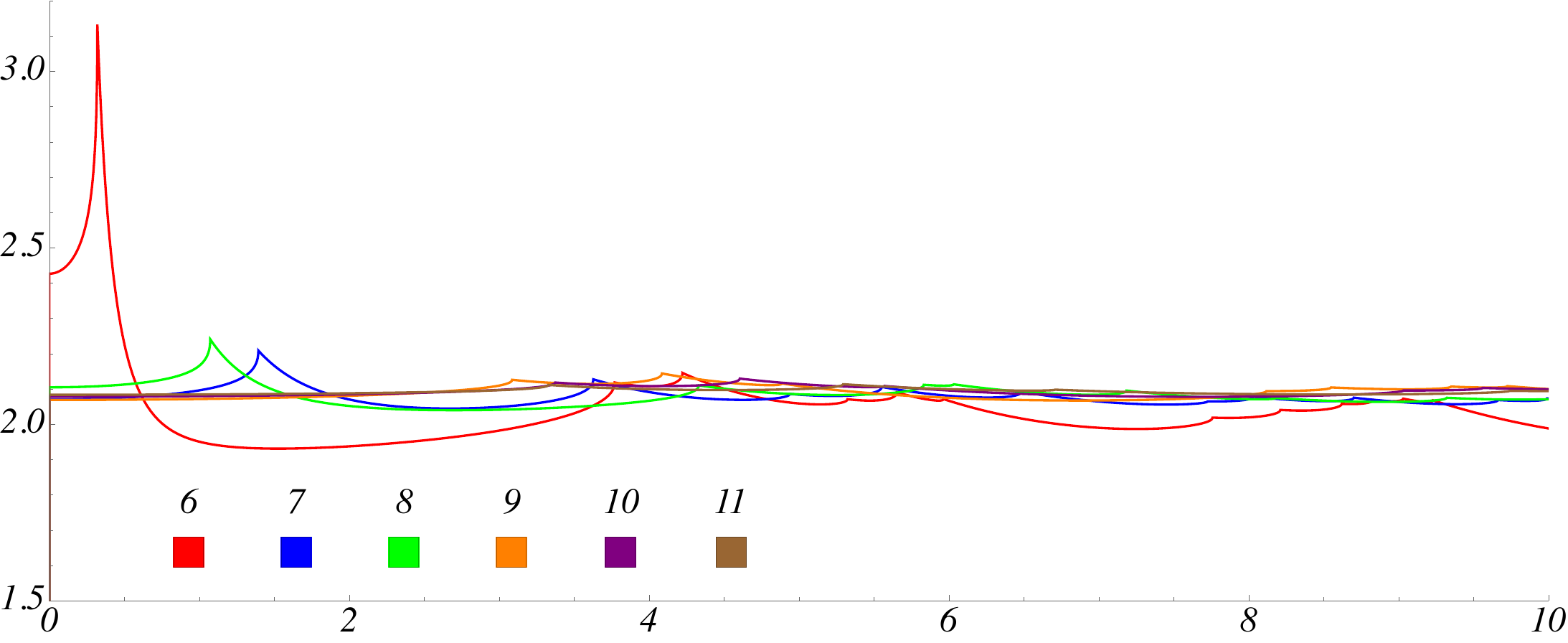}
\end{center}\vspace*{-0.7cm}
\begin{center}
\addtocounter{fig}{1} {\small Figure \arabic{fig}~: Graph of
  $\Theta_N$ for $N=10^m$ where $m$ is $6$ (red), $7$ (blue) up to
  $11$ (brown).}
\end{center}

Recall that the function $g$ defined in Equation \eqref{eq:defig}
satisfies $g(u)=\frac{2}{3} +\bigO(u^2)$ near $u=0$. Hence for all
$t\in [0,A]$ and $p\in \OOO_K\ssm\{0\}$, we have
\[
g\Big(\frac{t\,N}{2\,\phi(N)\,|p|} \Big) =
\frac{2}{3} +\bigO\Big(\frac{A^2\,\eta_N^{\;2}}{|p|^2}\Big)\,.
\]
By the right part of Equation \eqref{eq:controlcppvp} and by the fact
that the series $\sum_{p\in\OOO_K\smallsetminus\{0\}}\frac{1}{|p|^3}$
converges, we have
\[\sum_{p\in\OOO_K:\;0<|p|\leq N^\alpha}\frac{1}{c'_p\,|v_p|\,|p|^2}
=\bigO\Big(\sum_{p\in\OOO_K:\;0<|p|\leq N^\alpha}\frac{1}{|p|^3}\Big)
=\bigO(1)\,.
\]
Hence with $S_{N,0}$ the sum defined in Equation \eqref {eq:defiSN}
for $k=0$, for every $t\in [0,A]$, Equation \eqref{eq:defiThetaN}
gives
\begin{align*}
  \Theta_N(t)&=\frac{2\,N^3}{3\,\psi(N)\,\phi(N)}
  \sum_{\substack{p\in \OOO_K\\0<|p|\leq N^\alpha}} \frac{1}{c'_p\,|v_p|}+
  \bigO\Big(\frac{A^2\,\eta_N^{\;2}\,N^3}{\psi(N)\,\phi(N)}
  \sum_{\substack{p\in \OOO_K\\0<|p|\leq N^\alpha}}
  \frac{1}{c'_p\,|v_p|\,|p|^2}\Big)\\&=
\frac{8\,\pi\,N^3\,S_{N,0}}{3\,|D_K|\,\psi(N)\,\phi(N)}+
\bigO\Big(\frac{A^2\,\eta_N^{\;3}\,N^2}{\psi(N)}\Big)\,.
\end{align*}
Since the support of $f$ is contained in $[0,A]$, it follows from
Equation \eqref{eq:caspasplusinfty} that
\begin{align}
  \R^+_{N}(f)&=\frac{N^3\,S_{N,0}}{\psi(N)\;\phi(N)}\int_{0}^{+\infty} f(t)
  \frac{8\pi}{3\,|D_K|}\;dt+\bigO\Big(\frac{A^3\,\eta_N^{\;3}\,N^2}{\psi(N)}\,
  \|f\|_\infty\Big)\nonumber
  \\ &+\bigO\Big(\frac{\phi(N)\,N^{4\alpha+2\ga}}{\psi(N)}\,
  \|f'\|_\infty\Big)+\bigO\Big(\frac{N^{2+2\alpha-2\ga}}{\psi(N)}\,
  \|f\|_\infty\Big)\,.\label{eq:casplusinftyprov}
\end{align}
Since we have $\psi(N)=\frac{N^3\,S_{N,0}}{\phi(N)}$ in this Case
\hyperlink{cas5}{5}, since $\eta_N=\frac{N}{\phi(N)}$ remains
bounded, and by Equation \eqref{eq:controlSN} with $k=0$, we have
\begin{equation}\label{eq:errorterm1cas5}
\frac{\eta_N^{\;3}\,N^2}{\psi(N)}
=\frac{N^5\,\phi(N)}{\phi(N)^3N^3\,S_{N,0}}=
\bigO\Big(\frac{N^2}{\phi(N)^2\,N^\alpha}\Big)=
\bigO\Big(\frac{N^3}{\phi(N)^3}\,\frac{\phi(N)}{N^{1+\alpha}}\Big)=
\bigO\Big(\frac{\phi(N)}{N^{1+\frac{\alpha}{2}}}\Big)\,.
\end{equation}
Similarly, since $\ga=\frac{1-4\alpha}{2}$, ${\displaystyle
  \lim_{N\ra+\infty}}\;\frac{\phi(N)}{N^{1+\frac{\alpha}{2}}}=0$ and
$\alpha\leq\frac{2}{11}$ in this Case \hyperlink{cas5}{5}, we
have
\begin{equation}\label{eq:errorterm2cas5}
\frac{\phi(N)\,N^{4\alpha+2\ga}}{\psi(N)}
=\frac{\phi(N)^2\,N}{N^3\,S_{N,0}}=
\bigO\Big(\frac{\phi(N)^2}{N^{2+\alpha}}\Big)=
\bigO\Big(\frac{\phi(N)}{N^{1+\frac{\alpha}{2}}}\Big)^2=
\bigO\Big(\frac{\phi(N)}{N^{1+\frac{\alpha}{2}}}\Big)
\end{equation}
\begin{equation}\label{eq:errorterm4cas5}
\text{and}\qquad \frac{N^{2+2\alpha-2\ga}}{\psi(N)}
=\frac{\phi(N)\,N^{1+6\alpha}}{N^3\,S_{N,0}}=
\bigO\Big(\frac{\phi(N)}{N^{2-5\alpha}}\Big)=
\bigO\Big(\frac{\phi(N)}{N^{1+\frac{\alpha}{2}}}\Big)\,.
\end{equation}
Therefore, putting together Equations \eqref{eq:casplusinftyprov} to
\eqref{eq:errorterm4cas5}, we have
\begin{align}
  &\R^+_{N}(f) =\int_{0}^{+\infty} f(t)
  \frac{8\pi}{3\,|D_K|}\;dt+\bigO\Big(\frac{\phi(N)}{N^{1+\frac{\alpha}{2}}}
  \,(\|f'\|_\infty+A^3\|f\|_\infty)\Big)\,.\label{eq:casplusinfty}
\end{align}
By symmetry, for every $f\in C_c(\RR)$, this proves the last case of
Theorem \ref{theo:sumsquaremain}.

When $\phi:N\ra N^\beta$ is a power function, the assumptions of Case
\hyperlink{cas5}{5} on $\phi$ are satisfied if and only if $\beta\in
[1,1+ \frac{\alpha}{2}[\,$. If furthermore $D_K\equiv 0\!\!\mod 4$,
then by Equation \eqref{eq:defiSNcasdiscnul} with $k=0$, as
$N\ra+\infty$, we have $\psi(N)=\frac{N^3\,S_{N,0}}{\phi(N)}\sim
N^{3+\alpha-\beta}$. This proves the last case of Theorem
\ref{theo:mainintro} when $\beta\in[1,1+\frac{\alpha}{2}[\,$.

\medskip
\noindent{\bf Case \hypertarget{cas4}{4}. } Let us assume that we have
${\displaystyle\lim_{N\ra+\infty}}\frac{N^{1-\frac{\alpha}{2}}}{\phi(N)}
=0$, that $\lambda''_N={\displaystyle\lim_{N\ra+\infty}}
\frac{\phi(N)}{N}=0$ and that $\alpha\leq\frac{2}{11}$.

Note that the second assumption implies that $\lambda'''_N=
{\displaystyle \lim_{N\ra+\infty}}\;\frac{\phi(N)}
{N^{1+\frac{\alpha}{2}}} =0$. As in Case \hyperlink{cas5}{5}, we take
$\psi:N\mapsto \frac{N^3\,S_{N,0}}{\phi(N)}$ and
$\ga=\frac{1-4\alpha}{2}$. The assumptions \eqref{eq:conditiongaeps}
and \eqref{eq:scalingrange} are satisfied (since $\frac{\phi(N)}
{N^{2-2\alpha-2\ga}}= \frac{\phi(N)}{N}\,N^{-2\alpha}$ tends to $0$ as
$N\ra+\infty$ under the assumptions of Case \hyperlink{cas4}{4}),
hence we may apply the results of Section \ref{sec:linear}. We are
going to prove that also in this case, the function $\Theta_N$ defined
in Equation \eqref{eq:defiThetaN} converges uniformly on compact
subsets of $[0,+\infty[$ to the constant function $\frac{8\pi}
{3|D_K|}$. Below are the graphs of $\Theta_N$ for various $N$, in the
Gaussian case $K=\QQ[i]$, with $\alpha=0.15$ and $\beta =0.95\in
\;]1-\frac{\alpha}{2},1[\,$.
%, and in black the horizontal line with height
%$\frac{8\pi}{3|D_K|} \simeq 2.094$.
\begin{center}
    \includegraphics[width=15cm]{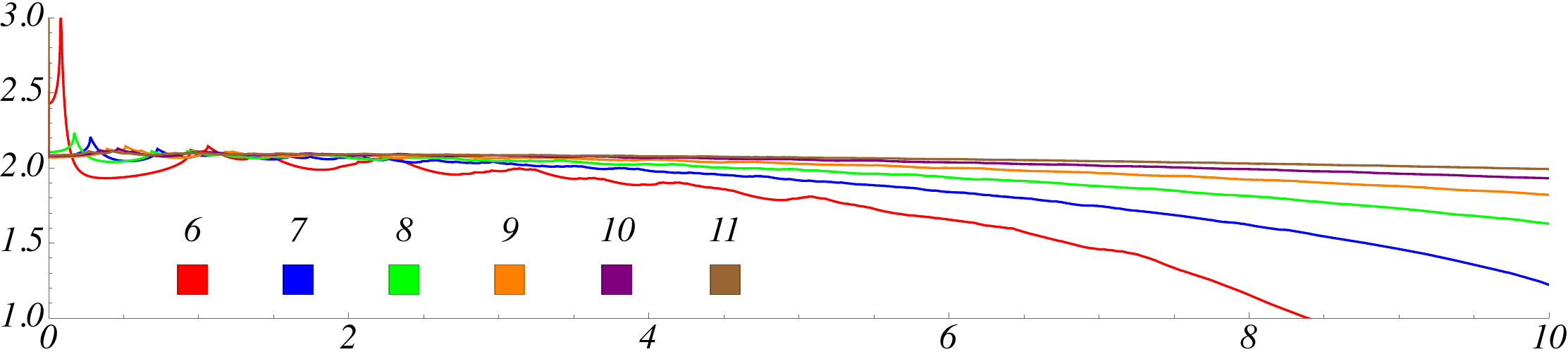}
\end{center}\vspace*{-0.7cm}
\begin{center}
\addtocounter{fig}{1} {\small Figure \arabic{fig}~: Graph of
  $\Theta_N$ for $N=10^m$ where $m$ is $6$ (red), $7$ (blue) up to
  $11$ (brown).}
\end{center}

The proof of Case \hyperlink{cas4}{4} is almost the same one as the
one of Case \hyperlink{cas5}{5}. Note that the quantity
$\eta_N=\frac{N} {\phi(N)}$ is now no longer bounded. Equation
\eqref{eq:errorterm1cas5} needs to be replaced by
\begin{equation}\label{eq:errorterm1cas4}
\frac{\eta_N^{\;3}\,N^2}{\psi(N)}
=\frac{N^5\,\phi(N)}{\phi(N)^3N^3\,S_{N,0}}=
\bigO\Big(\frac{N^2}{\phi(N)^2\,N^\alpha}\Big)=
\bigO\Big(\frac{N^{1-\frac{\alpha}{2}}}{\phi(N)}\Big)^2=
\bigO\Big(\frac{N^{1-\frac{\alpha}{2}}}{\phi(N)}\Big)\,.
\end{equation}
This is possible by the first assumption of Case \hyperlink{cas4}{4}
requiring that ${\displaystyle\lim_{N\ra+\infty}}
\frac{N^{1-\frac{\alpha}{2}}}{\phi(N)}=0$. Under the second assumption
of Case \hyperlink{cas4}{4}, we have
\[
\frac{\phi(N)}{N^{1+\frac{\alpha}{2}}}=\frac{N^{1-\frac{\alpha}{2}}}{\phi(N)}
\Big(\frac{\phi(N)}{N}\Big)^2=
\bigO\Big(\frac{N^{1-\frac{\alpha}{2}}}{\phi(N)}\Big)\,.
\]
Hence Equation \eqref{eq:errorterm2cas5} (which is still valid) gives
\begin{equation}\label{eq:errorterm2cas4}
  \frac{\phi(N)\,N^{4\alpha+2\ga}}{\psi(N)}=
\bigO\Big(\frac{N^{1-\frac{\alpha}{2}}}{\phi(N)}\Big)\,.
\end{equation}
The limit distribution is still Poissonian with the same constant
asymptotic pair correlation density given by $\frac{8\,\pi}
{3\,|D_K|}$.  Only the error term that involves the norm
$\|f\|_\infty$ is weakened, corresponding to the new Equations
\eqref{eq:errorterm1cas4} and \eqref{eq:errorterm2cas4}, so that
Equation \eqref{eq:casplusinfty} becomes
\[
  \R^+_{N}(f) =\int_{0}^{+\infty} f(t)\frac{8\pi}{3\,|D_K|}\;dt
  +\bigO\Big(\frac{\phi(N)}{N^{1+\frac{\alpha}{2}}}\,\|f'\|_\infty
  +\frac{N^{1-\frac{\alpha}{2}}}{\phi(N)}\,A^3\,\|f\|_\infty\Big)\,.
\]
By symmetry, for every $f\in C_c(\RR)$, this proves the penultimate
case of Theorem \ref{theo:sumsquaremain}. 

When $\phi:N\ra N^\beta$ is a power function, the assumptions of Case
\hyperlink{cas4}{4} on $\phi$ are satisfied if and only if
$\beta\in\;]1-\frac{\alpha}{2}, 1[\,$. If furthermore $D_K\equiv
0\!\!\mod 4$, then by Equation \eqref{eq:defiSNcasdiscnul} with $k=0$,
as $N\ra+\infty$, we have $\psi(N)=\frac{N^3\,S_{N,0}}{\phi(N)}\sim
N^{3+\alpha-\beta}$. This proves the last case of Theorem
\ref{theo:mainintro} when $\beta\in\; ]1-\frac{\alpha}{2},1[\,$.

This concludes the proof of Theorem \ref{theo:sumsquaremain}, and
along the way, we proved Theorem \ref{theo:mainintro} in the
Introduction.  \cqfd

\section{Experiments for bigger scalings}
\label{rem:casBetageq2}

The method used in order to prove the Poissonian asymptotic behaviour
of the pair correlations in Theorem \ref{theo:sumsquaremain} does not
work when $\lambda'''_\phi>0$ (with the notation of the beginning of
Section \ref{sec:mainproof}). Recall that this means that $\beta>1+
\frac{\alpha}{2}$ for power scalings $\phi:N\mapsto N^\beta$. We only
proved in Theorem \ref{theo:mainintro} that we have a Poissonian
asymptotic behaviour when $\beta\in\;]1-\alpha,1+\frac{\alpha}{2}[\,$.

Figure \hyperlink{figlastsection}{\addtocounter{fig}{1}\arabic{fig}%
  \addtocounter{fig}{-1}} below shows the empirical pair correlation
measures $\R_{2000}^{\alpha,\beta}$ (defined just before Theorem
\ref{theo:mainintro} in the Introduction) in the Gaussian case
$K=\QQ(i)$ with $\alpha=0.15$ for values
$\beta\in\{1.3,1.5,1.7,1.9,2\}$ outside the range of Theorem
\ref{theo:mainintro}. Again the distributions are renormalized by
$\psi(N)= N^{3+\alpha-\beta}$. These graphs indicate that, for some
parameters in the range from $1+\frac\alpha 2$ to $2$, possibly
including the interval $[1.7,2]$, the asymptotic pair correlation
densities should exist and should exhibit a {\it strong level
  repulsion}, that is, they should vanish on an open interval around
$0$.

\begin{center}
\begin{overpic}[width=14cm,tics=10]{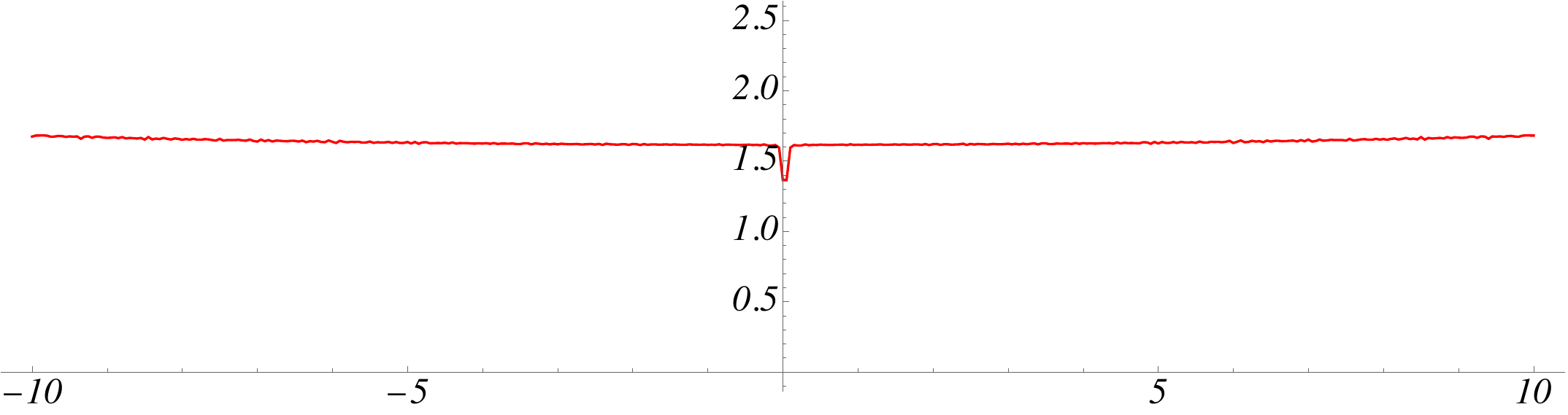}
 \put (85,10) {$\beta=1.3$}
\end{overpic}
\begin{overpic}[width=14cm,tics=10]{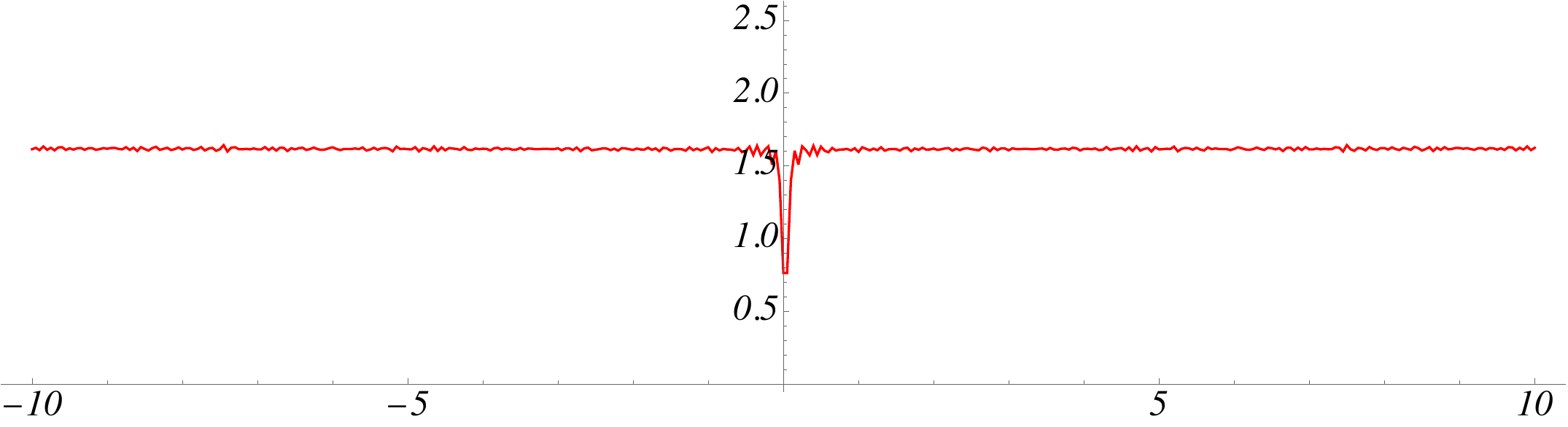}
 \put (85,10) {$\beta=1.5$}
\end{overpic}
\begin{overpic}[width=14cm,tics=10]{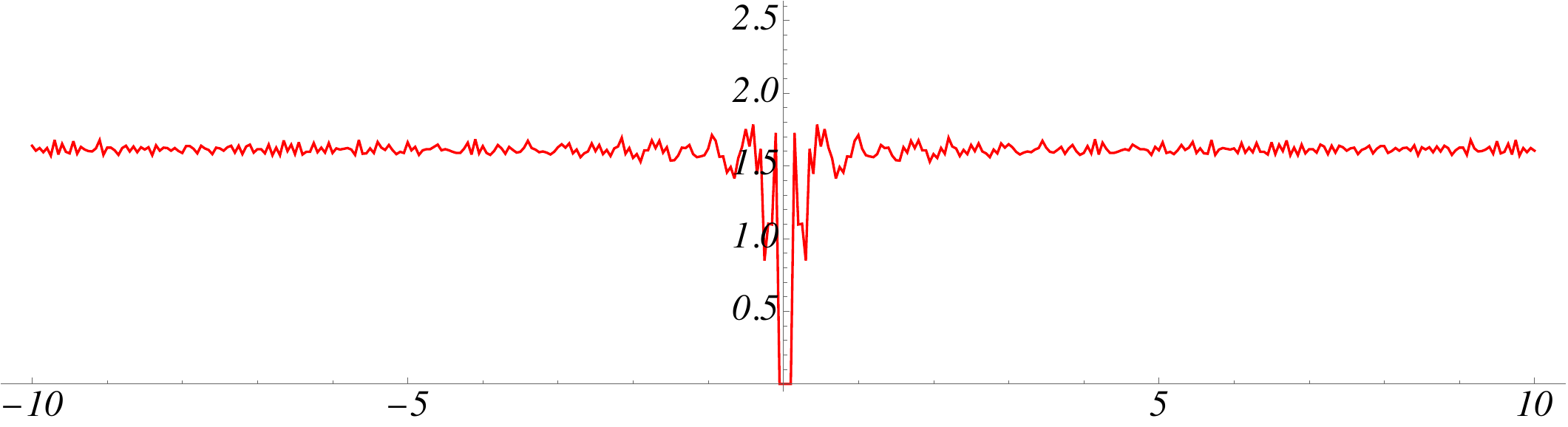}
 \put (85,10) {$\beta=1.7$}
\end{overpic}
\begin{overpic}[width=14cm,tics=10]{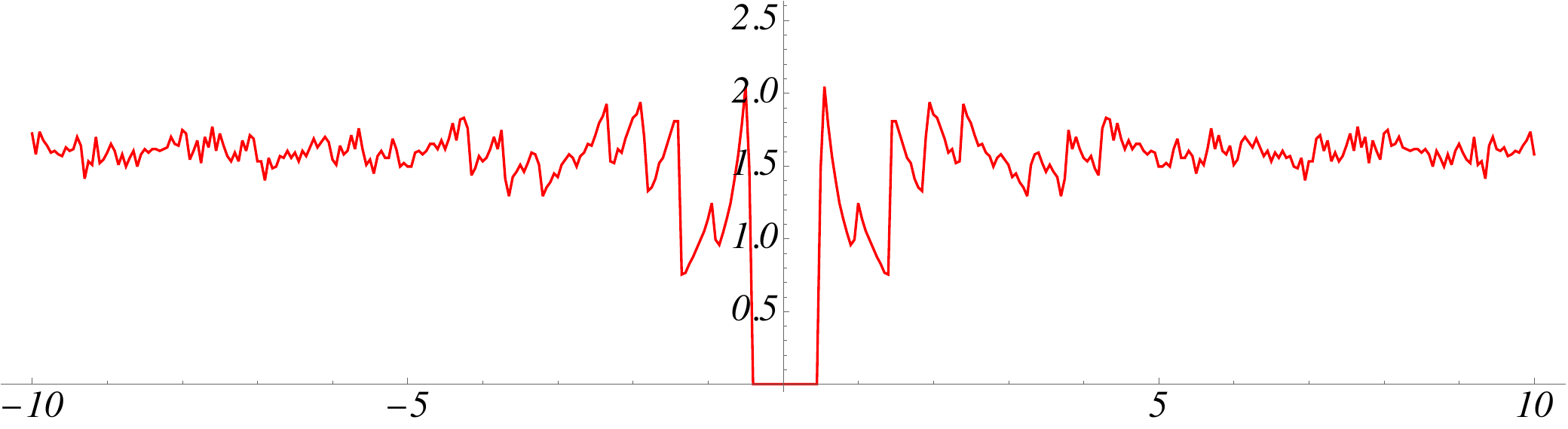}
 \put (85,10) {$\beta=1.9$}
\end{overpic}
\begin{overpic}[width=14cm,tics=10]{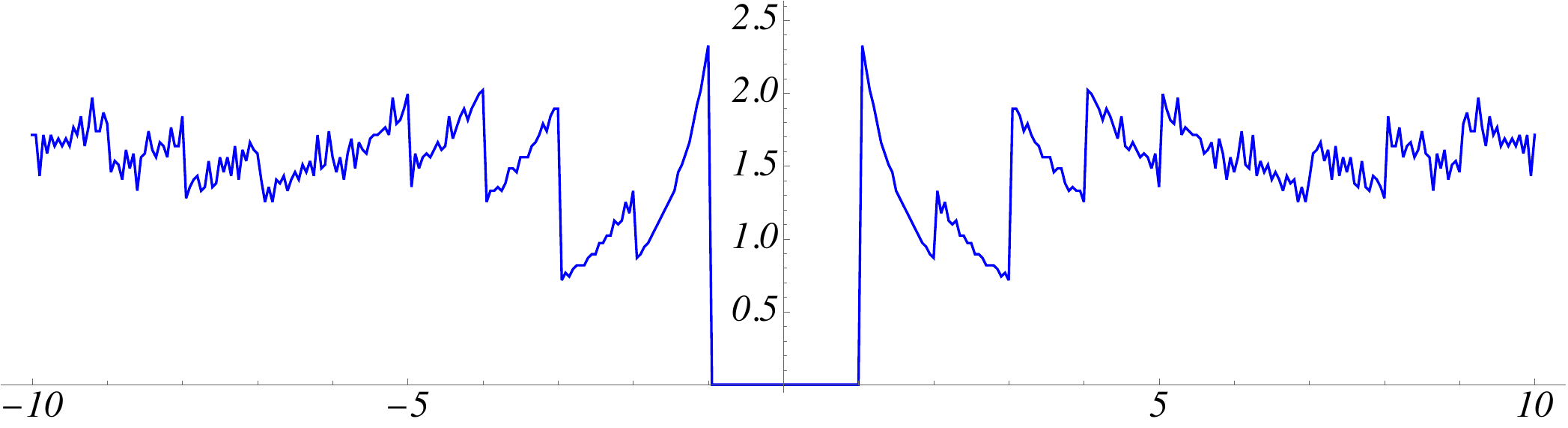}
 \put (85,8) {$\beta=2$}
\end{overpic}
     \end{center}\vspace*{-0.7cm}
\hypertarget{figlastsection}{\begin{center}
\addtocounter{fig}{1}
{\small Figure \arabic{fig}~: The empirical pair correlation
distributions $\R_{2000}^{0.15,\beta}$ with $\beta\in\{1.3,1.5,1.7,1.9,2\}$.}
\end{center}}

Note that the minimal difference of two sums of squares $a<b\le N^2$
is $1$, hence we have
\[
N^2(\ln b-\ln a)=
N^2\ln\Big(1+\frac{b-a}a\Big)\ge N^2\ln\Big(1+\frac 1{N^2}\Big)
\to_{N\to+\infty} 1\,.
\]
Thus, if a nonzero asymptotic pair correlation density exists for the
quadratic scaling (and we conjecture that it does), it will exhibit a
strong level repulsion, since it will vanish on $]-1,1[\,$. An analogous
computation indicates that there should be a total loss of mass for larger
power scalings $\phi:N\mapsto N^\beta$ with $\beta>2$.

%{\small \bibliography{viitteet} }

\bigskip
{\small
\noindent \begin{tabular}{l} 
Department of Mathematics and Statistics, P.O. Box 35\\ 
40014 University of Jyv\"askyl\"a, FINLAND.\\
{\it e-mail: jouni.t.parkkonen@jyu.fi}
\end{tabular}

\medskip
\noindent \begin{tabular}{l}
Laboratoire de math\'ematique d'Orsay, UMR 8628 CNRS,\\
Universit\'e Paris-Saclay,\\
91405 ORSAY Cedex, FRANCE\\
{\it e-mail: frederic.paulin@universite-paris-saclay.fr}
\end{tabular}}

\end{document}